\documentclass[10pt]{amsart}

\usepackage{amsmath}
\usepackage{amssymb}
\usepackage{amsfonts}
\usepackage{amsthm}
\usepackage{enumerate}
\usepackage[all]{xy}
\usepackage[latin1]{inputenc}
\usepackage{mathdots}
\usepackage{graphicx}
\usepackage{latexsym}
\usepackage{mathrsfs}
\usepackage{textcomp}
\usepackage{color}
\input xy

\makeatletter

\newtheoremstyle{definition}
        {5pt}
        {3pt}
        {}
        {0pt}
        {\scshape}
        {.}
        {5pt}
        {\thmname{#1} \thmnumber{#2} \thmnote{[#3]}} 

\newtheoremstyle{theorems}
        {5pt}
        {3pt}
        {\itshape}
        {0pt}
        {\scshape}
        {.}
        {5pt}
        {\thmname{#1} \thmnumber{#2}\thmnote{[#3]}} 

\swapnumbers

\renewcommand\section{\@startsection{section}{1}{\z@}%
        {-3.5ex \@plus -1ex \@minus -.2ex}%
        {2.3ex \@plus .2ex}%
        {\centering\reset@font\scshape}}

\makeatother

\theoremstyle{theorems}
\newtheorem{Theo}{Theorem}[section]
\newtheorem{Prop}[Theo]{Proposition}
\newtheorem{Cor}[Theo]{Corollary}
\newtheorem{Lemma}[Theo]{Lemma}

\theoremstyle{definition}

\newcommand{\rep}{{\rm rep}}

\newcommand{\Hom}{{\rm Hom}}

\newcommand{\Ext}{{\rm Ext}}
\newcommand{\End}{{\rm End}}

\newcommand{\DTr}{{\rm D\hskip -1pt Tr}\hskip 0.6 pt}
\newcommand{\TrD}{{\rm Tr \hskip -0.3 pt D}\hskip 0.6 pt}

\def\id{\hbox{1\hskip -3pt {\sf I}}}

\begin{document}

\title{\sc  Almost split sequences and approximations \vspace{5pt}}

\author{\sc Shiping Liu, Puiman Ng, and Charles Paquette}

\date{}

\address{Shiping Liu, D\'epartement de math\'ematiques, Universit\'e de Sherbrooke, Sherbrooke, Qu\'ebec, Canada, J1K 2R1}
\email{shiping.liu@usherbrooke.ca}

\address{Puiman Ng, D\'epartement de math\'ematiques, Universit\'e de Sherbrooke, Sherbrooke, Qu\'ebec, Canada, J1K 2R1}
\email{pui.man.ng@usherbrooke.ca}

\address{Charles Paquette, Department of Mathematics and Statistics, University of New Brunswick , Fredericton, New Brunswick, Canada}
\email{charles.paquette@usherbrooke.ca}

\maketitle

\begin{abstract}

\vspace{-20pt}

Let $\mathcal A$ be an exact category, that is, an extension-closed full subcategory of an abelian
category. Firstly, we give some necessary and sufficient conditions for $\mathcal A$ to have almost split sequences. Then, we study when an almost split sequence in $\mathcal A$ induces an almost split sequence in an exact subcategory  $\mathcal C$ of $\mathcal A$. In case $\mathcal A$ has almost split sequences and $\mathcal C$ is Ext-finite Krull-Schmidt, this provides a necessary and sufficient condition for $\mathcal C$ to have almost split sequences. Finally, we show two applications of these results.

\end{abstract}

\maketitle

\bigskip
\bigskip

\centerline{\sc Introduction}

\bigskip\medskip

The Auslander-Reiten theory of almost split sequences has been playing a fundamental role in the representation theory of artin algebras with a great impact in other areas such as algebraic geometry and algebraic topology; see \cite{ARS, Aus, Jor}.
It is a long standing problem to determine which categories have almost split sequences. In the mo\-dule category over an artin algebra, the existence of almost split sequences is derived from the Auslander-Reiten duality. In a general Hom-finite Krull-Schmidt exact category, Gabriel and Roiter showed that the existence of the Auslander-Reiten duality is necessary for the existence of almost split sequences; see \cite{GaR}, which is later proved to be sufficient by Lenzing and Zuazua in case the cate\-gory is in addition Ext-finite; see \cite{LeZ}. On the other hand, it is natural to study when a subcategory of a category having almost split sequences has almost split sequences. A pioneering work in this direction by Auslander and Smal{\o} shows that functorially finite subcategories of a module category have almost split sequences; see \cite{AuS}. Rather recently, J{\o}rgensen considered the analogous problem for Hom-finite Krull-Schmidt triangulated categories, and proved that if the ambient category has almost split triangles, then the almost split triangles in a subcategory are linked to those in the ambient category by minimal approximations; see \cite{Jor2}. By passing through the homotopy category, this result is applied to obtain a necessary and sufficient condition for a subcategory of a module category to have almost split sequences; see \cite{Ng}.

\bigskip

In this paper, we shall deal with the above mentioned two problems in a setup as general as possible. First, working with an arbitrary exact category, we shall characterize an almost split sequence in terms of linear forms on the stable endomorphism algebras of its end terms. This yields necessary and sufficient conditions for an exact category to have an almost split sequence with two prescribed end terms. Specializing to Hom-finite exact categories, we recover the above-mentioned result by Gabriel-Roiter and Lenzing-Zuazua. Then, we investigate the relation between the almost split sequences in an exact category and those in its exact subcategories. The result says in particular that if the ambient category has almost split sequences, then the almost split sequences in a Hom-finite Krull-Schmidt exact subcategory are precisely the minimal projectively or injectively stable approximations of the almost split sequences in the ambient category. This is a strengthened analogous version, by means of a very different approach, of J{\o}rgensen's result stated in \cite{Jor}. Since our categories do not necessarily have projective or injective objects, one can not simply apply J{\o}rgensen's result in our situation as is done in \cite{Ng}. As an application, given any torsion theory in an exact category having almost split sequences, we show that the torsion subcategory has right almost split sequences and the torsion-free subcategory has left almost split sequences. Finally, we shall apply our results to study almost split sequences in the category of finitely presented representations of an infinite quiver.

\smallskip

{\center \section{Preliminaries}}

\smallskip

Throughout the paper, $R$ stands for a commutative ring, which is not necessarily artinian unless otherwise explicitly stated. An $R$-category is a category in which the morphism sets are $R$-modules and the composition of morphisms is $R$-bilinear. Let $\mathcal A$ be an additive $R$-category, which will be called {\it Hom-finite} if its morphism modules are all of finite $R$-length. An idempotent endomorphism $e: X\to X$ is said to {\it split} in $\mathcal{A}$ if there exist morphisms $p: X\to Y$ and $q: Y\to X$ such that $e=qp$ and $p\,q=\id_Y$. Moreover, an object $X$ is called
{\it strongly indecompo\-sable} if $\,\End_{\hspace{-1pt}\mathcal{A}\hspace{-1pt}}(X)$ is local, and {\it Krull-Schmidt} if it is a finite direct sum of strongly indecomposable objects. Now, $\mathcal A$ is called
{\it Krull-Schmidt} if every non-zero objects is Krull-Schmidt. It is shown, from a functorial point of view, by Gabriel and Roiter that $\mathcal{A}$ is Krull-Schmidt if and only if every non-zero object has a semiperfect endomorphism ring and all the idempotent endomorphisms split in $\mathcal A$; see \cite[(3.3)]{GaR}. Here, we present an alternative elementary proof using of the following probably well known observation.

\bigskip

\begin{Lemma}\label{sp-ring}

If $A$ is a semiperfect ring, then the complete sets of orthogonal pri\-mitive idempotents in $A$ are pairwise conjugate up to per\-mutation.

\end{Lemma}

\noindent{\it Proof.} Let $\{e_1, \ldots, e_n\}$ and $\{f_1, \ldots, f_m\}$ be complete sets of orthogonal primitive idempotents in a semiperfect ring $A$. Then $A=Ae_1\oplus \cdots \oplus Ae_n=Af_1\oplus\cdots\oplus Af_m,$ where $Ae_i$ and $Af_j$ are indecomposable. By Theorem 27.11 in \cite{AnF}, $n=m$ and there exists a permutation $\sigma$ such that $Ae_i\cong Af_{\sigma(i)}\vspace{1pt}$, for $i=1, \ldots, n$. Thus $e_i=(e_ib_if_{\sigma(i)})(f_{\sigma(i)}a_ie_i)\vspace{1.3pt}$ and $f_{\sigma(i)}=(f_{\sigma(i)}a_ie_i)(e_ib_if_{\sigma(i)})$, with $a_i, b_i\in A$. Setting $a=\sum_{i=1}^n(f_{\sigma(i)}a_ie_i),\vspace{1pt}$ we see that $a^{-1}=\sum_{i=1}^n (e_ib_if_{\sigma(i)})\vspace{1pt}$ and $f_{\sigma(i)}=ae_ia^{-1}$.
The proof of the lemma is completed.

\bigskip

The following statement implies in particular the above mentioned result of Gabriel and Roiter.

\bigskip

\begin{Prop}\label{KS-idemp} 
If $X$ is an object in $\mathcal{A}$,
then $X=X_1\oplus\cdots \oplus X_n$ with $X_i$ strongly indecomposable if and only if $\hspace{0.3pt}\End_{\hspace{-1pt}\mathcal{A}\hspace{-1.2pt}}(X)$ is semiperfect and all its idempotents split$\,;$ and in this case, each direct summand of $X$ is a direct sum of objects of a subfamily of $\{X_1, \ldots, X_n\}$, which is its unique
decomposition into a direct sum of indecomposable objects up to isomorphism and permutation.

\end{Prop}

\noindent{\it Proof.} Put $E=\End_{\hspace{-1pt}\mathcal{A}\hspace{-1pt}}(X)$. Suppose that $E$ is semiperfect. Then $E$ has a complete set $\{e_1, \ldots, e_n\}$ of ortho\-gonal primitive idempotents such that the $e_iEe_i$ are local; see \cite[(27.6)]{AnF}.
If the idempotents in $E$ split in $\mathcal{A}$, then $X=X_1\oplus \cdots \oplus X_n$ with $\End_{\hspace{-1pt}\mathcal{A}\hspace{-1pt}}(X_i)\cong e_iEe_i\vspace{1pt}$.

Suppose now that $X=X_1\oplus \cdots \oplus X_n$, with canonical injections $q_i: X_i\to X$, canonical projections $p_i: X\to X_i$, and local rings $\End_{\hspace{-1pt}\mathcal{A}\hspace{-1pt}}(X_i)$. Setting $e_i=q_ip_i$, we get a complete set $\{e_1, \ldots, e_n\}$ of orthogonal primitive idempotents in $E$ such that $e_i E e_i\cong \End_{\hspace{-1pt}\mathcal{A}\hspace{-1pt}}(X_i)\vspace{1pt}$. In particular, $E$ is semiperfect; see \cite[(27.6)]{AnF}. Let $f$ be a non-zero idempotent in $E$.  Since $E/{\rm rad} E$ is semi-simple, every non-zero idempotent in $E$ is a sum of orthogonal primitive idempotents. 
In particular, there exists a complete set $\{f_1, \ldots, f_n\}$ of orthogonal primitive idempotents in $E$ such that $f=f_1+\cdots+f_r$, with $0<r\le n\vspace{1pt}$. By Lemma \ref{sp-ring}, there exists a permutation $\sigma$ and an invertible $a\in E$ such that
$f_i=a\hspace{0.3pt}e_{\sigma(i)}a^{-1}$, $i=1, \ldots, n$.
Set $L=X_{\sigma(1)}\oplus \cdots \oplus X_{\sigma(r)}$, and
\vspace{-2pt}$$p=(p_{\sigma(1)}, \cdots, p_{\sigma(r)})^Ta^{-1}: X\to L, \, \mbox{ and } \, q=a(q_{\sigma(1)}, \ldots, q_{\sigma(r)}): L \to X. \vspace{-2pt}$$
Then $f=qp$ and $p\hspace{0.5pt}q=\id_L$. That is, $f$ splits.

Next, assume that $M$ is a non-zero direct summand of $X$ with a canonical injection $u: M\to X$ and a canonical projection $v: X\to M$. Set $f=uv$, a non-zero idempotent in $E$. As seen above, there exists a permutation $\sigma$ and morphisms $p: X\to X_{\sigma(1)}\oplus \cdots \oplus X_{\sigma(r)}$ and $q: X_{\sigma(1)}\oplus \cdots \oplus X_{\sigma(r)}\to X$ such that $f=qp$ and $p\hspace{0.5pt}q=\id_{X_{\sigma(1)}\oplus \cdots \oplus X_{\sigma(r)}}.$ This yields an isomorphism $p\hspace{0.3pt}u: M\to X_{\sigma(1)}\oplus \cdots \oplus X_{\sigma(r)}\vspace{1pt}$.

Finally, in order to show the uniqueness of the decomposition of $M$, we need only to consider the case where $M=X$. Suppose that $X=Y_1\oplus \cdots \oplus Y_m$, with canonical injections $u_i: Y_i\to X$, canonical projections $v_i: X\to Y_i$, and indecomposable objects $Y_i$. This yields a complete set $\{w_1, \ldots, w_m\}$ of orthogonal idempotents in $E$, where $w_i=u_iv_i$. Since the idempotents split in $\mathcal{A}$, the $w_i$ are primitive. By Lemma \ref{sp-ring}, $n=m$ and we may assume that there exists $b\in A$ such that $w_i=b\hspace{0.6pt}e_ib^{-1}$, for $i=1, \ldots, n$. Then $v_ib\hspace{0.6pt}q_i: X_i\to Y_i$ is an isomorphism, $i=1, \ldots, n$.
The proof of the proposition is completed.

\medskip

The following statement is an evident observation.

\smallskip

\begin{Lemma}\label{KS}

Let $\mathcal{I}$ be an ideal of $\mathcal A$, and let $X$ be an object in $\mathcal A$.

\begin{enumerate}[$(1)$]

\item If ${\rm End}_{\hspace{-1pt}\mathcal{A}\hspace{-1pt}}(X)$ is local, then $\mathcal{I}(X, X)=\End_{\hspace{-1pt}\mathcal{A}\hspace{-1pt}}(X)\vspace{1pt}$ or $\mathcal{I}(X, X)\subseteq {\rm rad}(\End_{\hspace{-1pt}\mathcal{A}\hspace{-1pt}}(X))\vspace{1.5pt}.$

\item If $\mathcal{A}$ is Krull-Schmidt, then so is ${\mathcal A} /{\mathcal I}$.

\end{enumerate}

\end{Lemma}


\medskip

To conclude this section, we recall some terminology which will be needed later. A morphism $f: X\to Y$ in $\mathcal{A}$ is {\it right minimal} if any endomorphism $g: X\to X$ such that $fg=f$ is an automorphism; and {\it left minimal} if any endomorphism $h: Y\to Y$ such that $hf=f$ is an automorphism. Let $\mathcal{C}$ be a full subcategory of $\mathcal{A}$, and let $X$ be an object in $\mathcal{A}$. A morphism $f: M\to X$ with $M\in \mathcal{C}$ is called a {\it right $\mathcal{C}$-approximation} of $X$ if $\Hom_{\hspace{-1pt}\mathcal{A}\hspace{-1pt}}(L, f): \Hom_{\hspace{-1pt}\mathcal{A}\hspace{-1pt}}(L, M)\to \Hom_{\hspace{-1pt}\mathcal{A}\hspace{-1pt}}(L, X)$ is surjective for any $L\in \mathcal{C}$; and a {\it minimal right $\mathcal{C}$-approximation} if, in addition, $f$ is right minimal. Dually, a morphism
$g: X\to N$ with $N\in \mathcal{C}$ is called a {\it left $\mathcal{C}$-approximation} of $X$ if $\Hom_{\hspace{-1pt}\mathcal{A}\hspace{-1pt}}(f, L): \Hom_{\hspace{-1pt}\mathcal{A}\hspace{-1pt}}(N, L)\to \Hom_{\hspace{-1pt}\mathcal{A}\hspace{-1pt}}(X, L)$ is surjective for any $L\in \mathcal{C}$; and a {\it minimal left $\mathcal{C}$-approximation} if, in addition, $f$ is left minimal. Moreover, one says that $\mathcal{C}$ is {\it contravariantly finite} in $\mathcal{A}$ if every object in $\mathcal{A}$ has a right $\mathcal{C}$-approximation, {\it covariantly finite} in $\mathcal{A}$ if every object in $\mathcal{A}$ has a left $\mathcal{C}$-approximation, and {\it functorially finite} if it is both contravariantly finite and covariantly finite in $\mathcal A$; see \cite{ASm}.

\smallskip

{\center \section{Existence of almost split sequences}}

\smallskip

For the rest of the paper, let $\mathcal{A}$ stand for an {\it exact} $R$-category, that is, an extension-closed full subcategory of an abelian $R$-category $\mathscr{A}$. Clearly, the idempotents in $\mathcal{A}$ split if and only if $\mathcal{A}$ is closed under direct summands in $\mathscr{A}$. The objective of this section is to study the existence of almost split sequences in $\mathcal A$.

\medskip

Let $\delta: \xymatrixcolsep{14pt}\xymatrix{0\ar[r] & X\ar[r]^f & Y \ar[r]^g& Z\ar[r]& 0}\vspace{0pt}$ be a short exact sequence in $\mathcal{A}$. We shall call $f$ a {\it proper monomorphism} and $g$ a {\it proper epimorphism} in $\mathcal A$. Given any morphisms $u: X\to M$ and $v: N\to Z$ in $\mathcal{A}$, since  $\mathcal{A}$ is extension-closed in $\mathscr{A}$, we have a pushout diagram \vspace{2pt}
$$\xymatrixcolsep{16pt}\xymatrixrowsep{16pt}\xymatrix{
\delta: & 0\ar[r] & X\ar[r]\ar[d]_u & Y \ar[d]\ar[r] & Z\ar[r]\ar@{=}[d]& 0\\
u\delta: & 0\ar[r] & M\ar[r] & L \ar[r] & Z\ar[r]& 0
} \vspace{2pt} $$ as well as a pullback diagram \vspace{5pt}
$$\xymatrixcolsep{16pt}\xymatrixrowsep{16pt}\xymatrix{
\delta v: & 0\ar[r] & X\ar[r]\ar@{=}[d] & E \ar[d]\ar[r] & N\ar[r]\ar[d]^v& 0\\
\delta: & 0\ar[r] & X\ar[r] & Y \ar[r] & Z\ar[r]& 0.
}$$
Thus, the equivalent classes of short exact sequences $\xymatrixcolsep{14pt}\xymatrix{0\ar[r] & X\ar[r] & E \ar[r] & Z\ar[r]& 0}$ in $\mathcal{A}$ form an abelian group $\Ext_{\hspace{-2pt}\mathcal{A}\hspace{-1pt}}^1(Z, X)$ under Baer sum, which is an ${\rm End}_{\hspace{-1pt}\mathcal{A}\hspace{-1pt}}(X)$-${\rm End}_{\hspace{-1pt}\mathcal{A}\hspace{-1pt}}(Z)$-bimodule under
the multiplications illustrated in the above diagrams. Since $\End_{\hspace{-1pt}\mathcal{A}\hspace{-1pt}}(X)$ is an $R$-algebra, $\Ext_{\hspace{-2pt}\mathcal{A}\hspace{-1pt}}^1(X, Z)$ is an $R$-module. We shall say that $\mathcal A$ is {\it Ext-finite} if $\Ext^1_{\hspace{-2pt}\mathcal{A}\hspace{-1pt}}(X, Y)$ is of finite $R$-length for all $X, Y\in \mathcal{A}$.

\bigskip

The stable categories of $\mathcal{A}$ introduced by Gabriel and Roiter are essential for our investigation; see \cite[(9.2)]{GaR}, and also \cite[(2.1)]{LeZ}. A morphism $u: X\to Y$ in $\mathcal A$ is called {\it injectively
trivial} if the map \vspace{2pt}
$$\Ext_{\hspace{-2pt}\mathcal{A}\hspace{-1pt}}^1(L, u): \, \Ext_{\hspace{-2pt}\mathcal{A}\hspace{-1pt}}^1(L, X)\to \Ext_{\hspace{-2pt}\mathcal{A}\hspace{-1pt}}^1(L, Y): \eta\mapsto u\eta$$ is zero for any $L\in \mathcal{A}$, or equivalently, $u$ factors through any proper monomorphism $v: X\to M$ in $\mathcal{A}$. Dually, $u: X\to Y$ is called {\it projectively trivial} if the map
\vspace{2pt} $$\Ext_{\hspace{-2pt}\mathcal{A}\hspace{-1pt}}^1(u, L): \, \Ext_{\hspace{-2pt}\mathcal{A}\hspace{-1pt}}^1(Y, L)\to \Ext_{\hspace{-2pt}\mathcal{A}\hspace{-1pt}}^1(X, L): \zeta\mapsto \zeta u$$ is zero for any $L\in \mathcal{A}$, or equivalently, $u$ factors through any proper epimorphism $w: N\to Y$ in $\mathcal{A}$. It is easy to verify that the injectively trivial morphisms in $\mathcal{A}$ form an ideal, denoted as $I_{_{\hspace{-2pt}\mathcal{A}}}$; and the projectively trivial morphisms form an ideal, as denoted as $P_{_{\hspace{-2pt}\mathcal{A}}}$. Now, the quotient category $\hspace{2pt}\overline{\hspace{-2pt} \mathcal A \hspace{-0.5pt}}\hspace{0.5pt}=\mathcal{A}/I_{_{\hspace{-2pt}\mathcal{A}}}$ is called
the {\it injectively stable category} of $\mathcal{A}$, while $\underline{\hspace{0.3pt}\mathcal A\hspace{-0.3pt}}\hspace{0.5pt}=\mathcal{A}/P_{_{\hspace{-2pt}\mathcal{A}}}$ is called the {\it projectively stable category}. In the sequel, for a morphism $u: X\to Y$ in $\mathcal A$, we shall denote by $\hspace{1pt}\overline{u\hspace{-0.5pt}}\hspace{0.8pt}$ and $\underline{u\hspace{-2pt}}\hspace{2pt}$ its images in ${\rm Hom}_{\hspace{1.5pt}\overline{\hspace{-1.5pt} \mathcal A \hspace{-0.8pt}}\hspace{1pt}}(X, Y)$ and ${\rm Hom}_{\hspace{0.5pt}\underline{\hspace{0.5pt}\mathcal A \hspace{-0.5pt}}\hspace{0.5pt}}(X, Y)$, respectively.
Finally, an object $X\in \mathcal{A}$ is called {\it Ext-injective} if every proper monomorphism $f: X\to Y$ is a section; and {\it Ext-projective} if every proper epimorphism $g: Y\to X$ is a retraction. It is easy to see that $X$ is Ext-injective if and only if $\id_X$ is injectively trivial, or equivalently, $X$ is zero in $\hspace{3pt}\overline{\hspace{-2pt} \mathcal A \hspace{-0.6pt}}\hspace{0.6pt}$. Dually,
$X$ is Ext-projective if and only if $\id_X$ is projectively trivial, or equivalently, $X$ is zero in $\hspace{0.5pt}\underline{\hspace{0.5pt}\mathcal A \hspace{-0.5pt}}\hspace{0.6pt}$.

\bigskip

Next, we recall from \cite{AuR} some terminology and facts for the Auslander-Reiten theory. Let $f: X\to Y$ be a morphism. One says that $f$ is {\it right almost split} if $f$ is not a retraction and every non-retraction morphism $g: M\to Y$ factors through $f$; and {\it minimal right almost split} if $f$ is right minimal and right almost split. In a dual manner, one defines $f$ to be {\it left almost split} and {\it minimal left almost split}. Note that if $f: X\to Y$ is left or right almost split, then $X$ or $Y$ is strongly indecomposable, respectively; see \cite{AuR}. A short exact sequence \vspace{-2pt}$$\delta: \xymatrixcolsep{16pt}\xymatrix{0\ar[r] & X\ar[r]^f & Y \ar[r]^g& Z\ar[r]& 0}\vspace{-1pt}$$ is called {\it almost split} if $f$ is minimal left almost split and $g$ is minimal right almost split; see \cite{AuR}. In this case, both the $\End_{\hspace{-1pt}\mathcal{A}\hspace{-1pt}}(X)$-socle and the $\End_{\hspace{-1pt}\mathcal{A}\hspace{-1pt}}(Z)$-socle of $\Ext_{\hspace{-2pt}\mathcal{A}\hspace{-1pt}}^1(Z, X)\vspace{1pt}$ are simple generated by $\delta$. Moreover, since $\delta$ is unique up to isomorphism for $X$ and for $Z$, we may write $X=\tau_{_{\hspace{-2pt}\mathcal{A}}\hspace{-1pt}}Z\vspace{1pt}$ and $Z=\tau^-_{_{\hspace{-2pt}\mathcal{A}}}\hspace{-1.5pt}X.$ We shall say that $\mathcal A$ has {\it right almost split sequences} if every indecomposable object is either Ext-projective or the ending term of an almost split sequence, $\mathcal A$ has {\it left almost split sequences} if every indecomposable object is either Ext-injective or the starting term of an almost split sequence, and $\mathcal A$ has {\it almost split sequences} if it has both left and right almost split sequences.
The following result of Auslander and Reiten is originally stated for abelian categories. However, the proof stated in \cite[(2.13),(2.14)]{AuR} works for exact categories.

\medskip

\begin{Lemma} \label{ARS-lem}

A short exact sequence $\xymatrixcolsep{14pt}\xymatrix{0\ar[r] & X\ar[r]^f & Y \ar[r]^g& Z\ar[r]& 0}$ in $\mathcal A$ is almost split if and only if $f$ is left almost split and $\End_{\hspace{-1pt}\mathcal{A}\hspace{-1pt}}(Z)$ is local$\,;$ if and only if $g$ is right almost split and $\End_{\hspace{-1pt}\mathcal{A}\hspace{-1pt}}(X)$ is local.

\end{Lemma}

\bigskip

From now on, fix an injective co-generator $I$ for the category ${\rm Mod}\hspace{0.3pt}R$ of all $R$-modules, which will be minimal if $R$ is artinian. Then we have an exact endofunctor $D=\Hom_R(-, \; I): {\rm Mod}\hspace{0.3pt}R \to {\rm Mod}\hspace{0.3pt}R.$  For $\,U, V\in {\rm Mod}\hspace{0.3pt}R$, an $R$-bilinear form $<\; , \; >: \, U\times V\to I$ is called {\it non-degenerate} provided that, for any non-zero element $u\in U$, there exists some $v\in V$ such that $<\hspace{-2pt}u, v\hspace{-2pt}>\ne 0$, and for any non-zero element $v\in V$, there exists some $u\in U$ such that $<\hspace{-2pt}u, v\hspace{-2pt}>\ne 0$. Observe that every $R$-linear form $\varphi: \Ext_{\hspace{-2pt}\mathcal{A}\hspace{-1pt}}^1(Z, X)\to I$ determines, for each $L\in \mathcal{A}$, two $R$-bilinear forms:
\vspace{2pt}
$$<\; ,\; >_\varphi\,: \; {\rm Hom}_{\hspace{1pt}\overline{\hspace{-1pt}\mathcal{A}\hspace{-1pt}}\hspace{1pt}}(L, X)\times \Ext_{\hspace{-2pt}\mathcal{A}\hspace{-1pt}}^1(Z, L)\to I:
(\bar f, \eta)\mapsto \varphi(f\eta),
\vspace{-2pt}$$ and
\vspace{-3pt} $$_\varphi \hspace{-2pt} <\; ,\;>: \; \Ext_{\hspace{-2pt}\mathcal{A}\hspace{-1pt}}^1(L, X) \times {\rm Hom}_{\hspace{0.5pt}\underline{\mathcal A \hspace{-1pt}}}\hspace{1pt}(Z, L) \to I: (\zeta, \,\underline{g\hspace{-2pt}}\hspace{2pt})\mapsto
\varphi(\zeta g).\vspace{5pt}$$

On the other hand, if $\delta$ is a non-zero extension in $\Ext_{\hspace{-2pt}\mathcal{A}\hspace{-1pt}}^1(Z, X)$, then there exists always an $R$-linear form $\varphi: \Ext_{\hspace{-2pt}\mathcal{A}\hspace{-1pt}}^1(Z, X)\to I$ such that $\varphi(\delta)\ne 0$. We are now ready to state the following result of Gabriel and Roiter, which is implicitly stated in \cite[(9.3)]{GaR}; see also \cite[(3.1)]{LeZ}. We include their proof for the reader's convenience,

\bigskip

\begin{Prop}\label{bilinearform} Let $\,\delta: \xymatrixcolsep{14pt}\xymatrix{0\ar[r] & X\ar[r] & Y \ar[r] & Z\ar[r]& 0}$ be an almost split sequence in $\mathcal A$, and let $\varphi: \Ext_{\hspace{-2pt}\mathcal{A}\hspace{-1pt}}^1(Z, X)\to I$ be an $R$-linear form. If $\varphi(\delta)\ne 0$,  then the $R$-bilinear forms
\vspace{-5pt}
$$<\; ,\; >_\varphi\,: \; {\rm Hom}_{\hspace{1pt}\overline{\hspace{-2pt}\mathcal{A}\hspace{-0.6pt}}\hspace{0.6pt}}(L, X)\times \Ext_{\hspace{-2pt}\mathcal{A}\hspace{-1pt}}^1(Z, L)\to I
\vspace{-5pt}$$
and
\vspace{-1pt} $$_\varphi\hspace{-2pt} <\; ,\;>: \; \Ext_{\hspace{-2pt}\mathcal{A}\hspace{-1pt}}^1(L, X) \times {\rm Hom}_{\underline{\hspace{1pt}\mathcal{A}\hspace{-1pt}}}\hspace{1pt}(Z, L) \to I \vspace{2pt}$$
are both non-degenerate, for every $L\in \mathcal A$.

\end{Prop}

\noindent {\it Proof.} Suppose that $\varphi(\delta)\ne 0$ and $L\in \mathcal{A}.$ We shall prove only that $_\varphi\hspace{-2pt} <\; ,\;>$ is non-degenerate. Let $g: Z\to L$ be a morphism in $\mathcal A$ which is not projectively trivial. Then $\mathcal A$ admits a pullback diagram
\vspace{-0pt}
$$\xymatrixcolsep{16pt}\xymatrixrowsep{14pt}\xymatrix{\eta: & 0\ar[r] & M\ar[r]\ar@{=}[d] & N \ar[d]\ar[r] & Z\ar[r]\ar[d]^g& 0\\
\zeta: & 0\ar[r] & M\ar[r] & E \ar[r] & L\ar[r]& 0}$$ with non-split rows. Since $\delta$ is almost split, there exists a pushout diagram
$$\xymatrixcolsep{16pt}\xymatrixrowsep{14pt}\xymatrix{\eta: & 0\ar[r] & M\ar[r]\ar[d]^h & N \ar[d]\ar[r] & Z\ar[r]\ar@{=}[d]& 0\\
\delta: & 0\ar[r] & X\ar[r] & Y \ar[r] & Z\ar[r]& 0}$$ in $\mathcal A$. This yields $\delta=h \eta =h(\zeta g)=(h\zeta) g$. As a consequence, $h\zeta\in \Ext_{\hspace{-2pt}\mathcal{A}\hspace{-1pt}}^1(L, X)$ is such that $_\varphi\hspace{-2pt}<\hspace{-2pt} h\zeta, g\hspace{-2pt}>=\varphi((h\zeta)g)=\varphi(\delta)\ne 0\vspace{1pt}$.
On the other hand, consider a non-split short exact sequence $\zeta: \,\xymatrixcolsep{15pt}\xymatrix{0\ar[r] & X\ar[r] & E \ar[r] & L\ar[r]& 0}$ in $\mathcal A$. Since $\delta$ is almost split, there exists a pullback diagram
$$\xymatrixcolsep{16pt}\xymatrixrowsep{14pt}\xymatrix{\delta: & 0\ar[r] & X\ar@{=}[d] \ar[r]& Y \ar[d]\ar[r] & Z\ar[r]\ar[d]^g& 0\\
\zeta: & 0\ar[r] & X\ar[r] & E \ar[r] & L\ar[r]& 0}$$ in $\mathcal A$. Thus $\underline{\hspace{0.5pt}g\hspace{-2pt}}\hspace{1pt} \in {\rm Hom\hspace{1pt}}_{\underline{\mathcal A \hspace{-1pt}}}\hspace{1pt}(Z, L)\vspace{1pt}$ is such that $_\varphi\hspace{-2pt}<\hspace{-2pt}\zeta, \, \underline{\hspace{0.5pt}g\hspace{-2pt}}\hspace{0pt}>=\varphi(\zeta g)=\varphi(\delta)\ne 0.$ The proof of the proposition is completed.

\bigskip

Our first result will be a characterization of an almost split sequence.
We need to introduce some terminology.
Let $F: \mathcal{A}\to {\rm Mod}\hspace{0.3pt}R$ and $G: \mathcal{A}\to {\rm Mod}\hspace{0.3pt}R$ be covariant or contravariant $R$-linear functors. A {\it functorial monomorphism} $\phi: F\to G$ is a natural transformation with $\phi_X: F(X)\to G(X)$ being injective for all $X\in \mathcal{A}$. Moreover, if $U\in {\rm Mod}\hspace{0.3pt}R$, then an $R$-linear form $\phi: U\to I$ is called {\it almost vanishing} if it vanishes on ${\rm rad}\hspace{0.3pt}U$ but not on $U$.

\bigskip

\begin{Theo} \label{AR-sequence} Let $\mathcal A$ be an exact $R$-category, which has a short exact sequence $\,\delta:\xymatrixcolsep{14pt}\xymatrix{0\ar[r] & X\ar[r] & Y \ar[r] & Z\ar[r]& 0}$ with
$X, Z$ strongly indecomposable. The following statements are equivalent.

\vspace{-1pt}

\begin{enumerate}[$(1)$]

\item The sequence $\delta$ is an almost split sequence in $\mathcal{A}\vspace{0.5pt}$.

\item There exists a functorial monomorphism $\phi: \Ext_{\hspace{-2pt}\mathcal{A}\hspace{-1pt}}^1(Z, \, -)\to D\Hom_{\hspace{1pt}\overline{\hspace{-2pt}\mathcal{A}\hspace{-0.6pt}}} \hspace{0.6pt}(-, \, X)$ such that $\phi_X(\delta)$ is almost vanishing on $\End_{\hspace{0.5pt}\overline{\hspace{-2pt} \mathcal{A} \hspace{-0.6pt}}\hspace{0.8pt}}(X)$.

\vspace{2pt}

\item There exists a functorial monomorphism $\psi: \Ext_{\hspace{-2pt}\mathcal{A}\hspace{-1pt}}^1(-,\, X) \to  D\Hom_{\hspace{2pt}\hspace{-2pt}\underline{\mathcal{A}\hspace{-0.6pt}}\hspace{0.6pt}}(Z, \, -)$ such that $\psi_Z(\delta)$ is almost vanishing on $\End_{\hspace{-1pt}\underline{\hspace{0.8pt} \mathcal{A} \hspace{-0.6pt}}\hspace{0pt}}(Z)$.

\end{enumerate}
\end{Theo}

\noindent{\it Proof.} We shall prove only the equivalence of Statements (1) and (2). Since $\delta\ne 0$ in each of the statements, we may assume that $X$ is not Ext-injective, that is, $X$ is non-zero in $\hspace{1pt}\overline{\hspace{-2pt} \mathcal{A} \hspace{-0.5pt}}\hspace{0.6pt}.$ By Lemma \ref{KS}(1), ${\rm rad}(\End_{\hspace{0.5pt}\overline{\hspace{-2pt}\mathcal{A} \hspace{-1pt}}\hspace{1pt}}(X))={\rm rad}(\End_{\hspace{-1pt}\mathcal{A} \hspace{-1pt}}(X))/I_{\hspace{-1pt}\mathcal{A}\hspace{-1pt}}(X, X).\vspace{1pt}$

Assume first that $\delta$ is an almost split sequence. In particular, there exists an $R$-linear form $\varphi: {\rm Ext}^1(Z, X)\to I$ such that $\varphi(\delta)\ne 0$. Fix $L\in \mathcal{A}.$ By Proposition \ref{bilinearform}, we have a non-degenerate bilinear form
$$
<\; ,\; >_\varphi\,: \; {\rm Hom}_{\hspace{1.5pt}\overline{\hspace{-1.5pt} \mathcal{A} \hspace{-1pt}}\hspace{1.5pt}}(L, X)\times \Ext_{\hspace{-2pt}\mathcal{A}\hspace{-1pt}}^1(Z, L)\to I:
(\bar f, \eta)\mapsto \varphi(f\eta).
$$
This induces an $R$-linear monomorphism

\vspace{-6pt}

$$\phi_L: \Ext_{\hspace{-2pt}\mathcal{A}\hspace{-1pt}}^1(Z, \, L)\to D\Hom_{\hspace{1pt}\overline{\hspace{-2pt}\mathcal{A}\hspace{-0.6pt}}\hspace{0.6pt}} (L, \, X): \eta\mapsto \,<\hspace{-3pt}-\, , \eta\hspace{-2pt}>_\varphi,\vspace{1pt}$$ which is clearly natural in $L$. Since $\phi_X$ is injective, $\phi_X(\delta)\ne 0$. If
$\bar{f}\in {\rm rad}(\End_{\hspace{0.5pt}\overline{\hspace{-2pt} \mathcal{A} \hspace{-1pt}}\hspace{1pt}}(X))\vspace{1pt}$, then $f\in{\rm rad}({\rm End}_{\hspace{-1.5pt}\mathcal{A}\hspace{-1pt}}(X)).$ Since $\delta$ is almost split, we have $f\delta=0$. As a consequence,
$\phi_X(\delta)(f)=<\hspace{-3pt}f, \delta\hspace{-2pt}>_\varphi=\varphi(f\delta)=0.$ Thus, $\phi_X(\delta)$ is almost vanishing on $\End_{\hspace{1pt}\overline{\hspace{-1pt} \mathcal{A} \hspace{-1pt}}\hspace{1pt}}(X).$

\vspace{2.5pt}

Conversely, let
$\phi: \Ext_{\hspace{-2pt}\mathcal{A}\hspace{-1pt}}^1(Z, \, -)\to D\Hom_{\hspace{1pt}\overline{\hspace{-2pt}\mathcal{A}\hspace{-0.6pt}}} \hspace{0.6pt}(-, \, X)\vspace{1pt}$ be a functorial monomorphism such that $\phi_X(\delta)$ is almost vanishing on $\End_{\hspace{0.5pt}\overline{\hspace{-2pt} \mathcal{A} \hspace{-1pt}}\hspace{1pt}}(X)\vspace{1pt}$. Then, $\delta\ne 0.$ Let $u: X\to L$ be a morphism in $\mathcal A$ which is not a section. For any morphism $v: L\to X$, we have $vu\in {\rm rad}(\End_{_{\hspace{-0.5pt}\mathcal A\hspace{-0.5pt}}}(X))$, and hence, $\bar{v}\bar{u}\in {\rm rad}(\End_{\hspace{0.5pt}\overline{\hspace{-2pt} \mathcal{A} \hspace{-1pt}}\hspace{1pt}}(X)).\vspace{1pt}$ Thus $\phi_X(\delta)(\bar v\bar u)=0$, that is, $(D\Hom_{\hspace{1pt}\overline{\hspace{-2pt}\mathcal A \hspace{-0.5pt}}}\hspace{0.5pt}(u, X) \circ \phi_X)(\delta)=0$. In view of the commutative diagram
\vspace{-2pt}
$$\xymatrix{\Ext_{\hspace{-2pt}\mathcal{A}\hspace{-1pt}}^1(Z, X)\ar[d]_{\phi_X} \ar[rr]^-{\Ext_{\hspace{-1pt}\mathcal{A}\hspace{-1pt}}^1(Z, u)} && \Ext_{\hspace{-2pt}\mathcal{A}\hspace{-1pt}}^1(Z, L)\ar[d]^{\phi_L}\\ D\Hom_{\hspace{0pt}\overline{\hspace{-2pt}\mathcal{A} \hspace{-0.6pt}}\hspace{0.6pt}}(X, X)\ar[rr]^-{D\Hom_{\hspace{1pt}\overline{\hspace{-1.5pt}\mathcal{A} \hspace{-0.6pt}}\hspace{0.6pt}}(u, X)}
&&  D\Hom_{\overline{\hspace{-2pt}\mathcal{A}\hspace{-0.6pt}}\hspace{0.6pt}}(L, X),} \vspace{-2pt}$$
we see that $(\phi_L\circ \Ext^1_{\mathcal A \hspace{-1pt}}\hspace{1pt}(Z, u))(\delta)=0$. Since $\phi_L$ is
injective, $u\delta=\Ext^1_{\mathcal A \hspace{-1pt}}\hspace{1pt}(Z, u)(\delta)=0$. That is, $u$ factors through the monomorphism $X\to Y$ in $\delta$. By Lemma \ref{ARS-lem}, $\delta$ is an almost split sequence. The proof of the theorem is completed.

\bigskip

If $X, Y\in \mathcal{A}$, then $D\Hom_{\hspace{1pt}\overline{\hspace{-2pt}\mathcal{A}\hspace{-0.6pt}}\hspace{0.6pt}}(X, Y)\vspace{1.5pt}$ is an $\End_{\hspace{-1pt}\mathcal{A}\hspace{-1pt}}(X)$-
$\hspace{-3pt}\End_{\hspace{-1pt}\mathcal{A}\hspace{-1pt}}(Y)$-bimodule with multiplications defined, for $f\in \End_{\hspace{-0.5pt}\mathcal{A}\hspace{-1pt}}(X)$, $\theta\in D\Hom_{\hspace{0.5pt}\overline{\hspace{-1.5pt}\mathcal{A}\hspace{-0.6pt}}\hspace{1pt}}(X, Y)$, $g\in \End_{\hspace{-0.5pt}\mathcal{A}\hspace{-1pt}}(Y)$, by
$$f\theta g: \Hom_{\hspace{0.5pt}\overline{\hspace{-1.5pt}\mathcal{A}\hspace{-0.6pt}}\hspace{1pt}}(X, Y)\to I: \overline{h}\mapsto \theta(\overline{ghf}).$$
Similarly, $D\Hom_{\hspace{0.3pt}\underline{\hspace{0.3pt}\mathcal{A}\hspace{-0.8pt}}}(X, Y)$ is an $\End_{\hspace{-1pt}\mathcal{A}\hspace{-1pt}}(X)$-$\hspace{-1pt}\End_{\hspace{-1pt}\mathcal{A}
\hspace{-1pt}}(Y)$-bimodule.

\medskip

\begin{Theo} \label{AR-exist-cor} Let $\mathcal{A}\vspace{1pt}$ be an exact $R$-category, and let  $X, Z$ be strongly indecomposable objects in $\mathcal A$. The following statements are equivalent.

\vspace{-0pt}

\begin{enumerate}[$(1)$]

\item There exists an almost split sequence $\xymatrixcolsep{14pt}\xymatrix{0\ar[r] & X\ar[r] & Y \ar[r] & Z\ar[r]& 0}$ in $\mathcal{A}$.

\vspace{1.5pt}

\item The $\End_{\hspace{-1pt}\mathcal{A}\hspace{-1pt}}(X)$-socle of $\Ext_{\hspace{-2pt}\mathcal{A}\hspace{-1pt}}^1(Z, X)\vspace{1pt}$ is non-zero and there exists a functorial monomor\-phism $\phi: \Ext_{\hspace{-2pt}\mathcal{A}\hspace{-1pt}}^1(Z, \, -)\to D\Hom_{\hspace{0.5pt}\overline{\hspace{-1.5pt}\mathcal{A}\hspace{-0.6pt}}\hspace{1pt}}(-, \, X)$.

\vspace{3pt}

\item The $\End_{\hspace{-1pt}\mathcal{A}\hspace{-1pt}}(Z)$-socle of $\Ext_{\hspace{-2pt}\mathcal{A}\hspace{-1pt}}^1(Z, X)\vspace{1pt}$ is non-zero and there exists a functorial monomorphism $\psi: \Ext_{\hspace{-2pt}\mathcal{A}\hspace{-1pt}}^1(-,\, X) \to  D\Hom_{\hspace{1pt}\underline{\mathcal A\hspace{-1pt}}\hspace{1pt}}(Z, \, -)$.

\end{enumerate}\end{Theo}

\noindent{\it Proof.} We shall prove only the equivalences of Statements (1) and (2). Assume first that $\mathcal A$ has an almost split sequence $\delta: \xymatrixcolsep{15pt}\xymatrix{0\ar[r] & X\ar[r] & Y \ar[r] & Z\ar[r]& 0.}$ By Theorem \ref{AR-sequence}, there exists a functorial monomor\-phism $\phi: \Ext_{\hspace{-2pt}\mathcal{A}\hspace{-1pt}}^1(Z, \, -)\to D\Hom_{\hspace{0.5pt}\overline{\hspace{-1.5pt}\mathcal{A}\hspace{-0.6pt}}\hspace{1pt}}(-, \, X)$. Being almost split, $\delta$ a non-zero element in the $\End_{\hspace{-1pt}\mathcal{A}\hspace{-1pt}}(X)$-socle of $\Ext_{\hspace{-2pt}\mathcal{A}\hspace{-1pt}}^1(Z, X)\vspace{1pt}$.

Conversely, suppose that $\delta: \xymatrixcolsep{15pt}\xymatrix{0\ar[r] & X\ar[r] & Y \ar[r] & Z\ar[r]& 0}$ is a non-zero extension lying in the $\End_{\hspace{-1pt}\mathcal{A}\hspace{-1pt}}(X)$-socle of $\Ext_{\hspace{-2pt}\mathcal{A}\hspace{-1pt}}^1(Z, X)\vspace{1pt}$. In particular,
$X$ is not Ext-injective and $Z$ is not Ext-projective. Let $\phi: \Ext_{\hspace{-2pt}\mathcal{A}\hspace{-1pt}}^1(Z, \, -)\to D\Hom_{\hspace{1pt}\overline{\hspace{-2pt}\mathcal{A}\hspace{-0.6pt}}\hspace{1pt}}(-, \, X)\vspace{1.5pt}$ be a functorial monomorphism. Since $\phi$ is natural, $\phi_X: \Ext_{\hspace{-2pt}\mathcal{A}\hspace{-1pt}}^1(Z, X)\to D\End_{\hspace{0.5pt}\overline{\hspace{-1.5pt}\mathcal{A}\hspace{-0.8pt}}\hspace{1pt}}(X)\vspace{1pt}$ is $\End_{\hspace{0pt}\mathcal{A}\hspace{-0.6pt}}(X)$-linear. Hence, $\theta=\phi_X(\delta)$ is a non-zero element of
$D\End_{\overline{\hspace{-1.5pt}\mathcal{A}\hspace{-0.4pt}}\hspace{0.6pt}}(X)\vspace{2pt}$, which is
annihilated by ${\rm rad}(\End_{\hspace{-1pt}\mathcal{A}\hspace{-1pt}}(X))$.
If $\bar{f}\in {\rm rad}(\End_{\hspace{0.5pt}\overline{\hspace{-2pt}\mathcal{A}\hspace{-1pt}}\hspace{1pt}}(X))\vspace{1pt}$, then $f\in {\rm rad}(\End_{\hspace{-1pt}\mathcal{A} \hspace{-1pt}}(X))$ by Lemma \ref{KS}(1), and hence,
$\theta(\bar{f})=(f\theta)(\overline{\id})=0$. That is, $\theta$ is almost vanishing. By Theorem \ref{AR-sequence}, $\delta$ is almost split in $\mathcal A$. The proof of the theorem is completed.

\bigskip

For the rest of this section, we shall concentrate on the case where $R$ is artinian and $I$ is the minimal injective co-generator for ${\rm Mod}\hspace{0.3pt}R$. In this case, we have a duality $D={\rm Hom}_R(-, I): {\rm mod}\hspace{0.3pt}R \to {\rm mod}\hspace{0.3pt}R,$ where ${\rm mod}\hspace{0.4pt}R$ denotes the category of finitely generated $R$-modules.

\medskip

\begin{Lemma}\label{duality}

Let $R$ be artinian, and let $\,<\; ,\;>: \,U\times V\to I$ be a non-degenerate $R$-bilinear form, where $U, V\in {\rm Mod}\hspace{0.3pt}R$. If $\,U$ or $V$ is of finite $R$-length, then we have two $R$-linear isomorphisms
$$V\to DU: v\mapsto\, <\hspace{-2pt}-\, , v\hspace{-2pt}> \; \mbox{ and } \; U\to DV: u\mapsto\, <\hspace{-2pt}u, \, -\hspace{-2pt}>.$$

\end{Lemma}

\noindent {\it Proof.} Since $<\; , \;>$ is non-degenerate, we have two $R$-linear monomorphisms$\,:$
\vspace{-0pt} $$\phi: V\to DU: v\mapsto\, <\hspace{-2pt}-\, , v\hspace{-2pt}> \; \mbox{ and } \; \psi: U\to DV: u\mapsto\, <\hspace{-2pt}u, \, -\hspace{-2pt}>. \vspace{-0pt}$$

Suppose that $V$ has finite $R$-length $\ell_R(V)$. Then $\ell_R(DV)=\ell_R(V)$. Since $\psi$ is injective, $\ell_R(U)\le \ell_R(V)$. On the other hand, since $\phi$ is injective, we have $\ell_R(V)\le \ell_R(DU)=\ell_R(U)\le \ell_R(V)$. This yields $\ell_R(V)=\ell_R(DU)$, and hence,
$\ell_R(U)=\ell_R(DV)$. As a consequence, $\phi$ and $\psi$ are isomorphisms. Similarly, the result holds true if $\ell_R(U)$ is finite. The proof of the lemma is completed.

\bigskip

\begin{Lemma} \label{AR-iso} Let $\mathcal{A}\vspace{1pt}$ be an exact $R$-category where $R$ is artinian, which has an almost split sequence $\xymatrixcolsep{14pt}\xymatrix{0\ar[r] & X\ar[r] & Y \ar[r] & Z\ar[r]& 0.}$

\begin{enumerate}[$(1)$]

\item If $L\in \mathcal A$, then
$\Ext_{\hspace{-2pt}\mathcal{A}\hspace{-1pt}}^1(Z,L)\vspace{1pt}$ is of finite $R$-length if and only if so is ${\rm Hom}_{\hspace{1pt}\overline{\hspace{-2pt}\mathcal{A}\hspace{-0.5pt}}\hspace{0.8pt}}(L,X);$ and in this case, $\Ext_{\hspace{-2pt}\mathcal{A}\hspace{-1pt}}^1(Z,\, L) \cong  \hspace{0.3pt} D\Hom_{\hspace{1pt}\overline{\hspace{-2pt}\mathcal{A}\hspace{-0.6pt}}\hspace{0.6pt}}(L, \,X)$.

\vspace{3pt}

\item If $L\in \mathcal A$, then $\Ext_{\hspace{-2pt}\mathcal{A}\hspace{-1pt}}^1(L,X)\vspace{1pt}$ is  of finite $R$-length if and only if so is $\Hom_{\hspace{0.3pt}\underline{\mathcal{A}\hspace{-1pt}}\hspace{1pt}}(Z,L);$ and in this case, $\Ext_{\hspace{-2pt}\mathcal{A}\hspace{-1pt}}^1(L,\, X) \cong D \Hom_{\hspace{0.3pt}\underline{\mathcal{A}\hspace{-1pt}}\hspace{1pt}}(Z,\, L).$

\end{enumerate}\end{Lemma}

\noindent{\it Proof.} We shall prove only Statement (1). For any $L\in \mathcal{A}$, by Proposition \ref{bilinearform}, there exists a non-degenerate $R$-bilinear form $<\; ,\; >: \; {\rm Hom}_{\hspace{1pt}\overline{\hspace{-2pt}\mathcal{A}\hspace{-0.6pt}}\hspace{0.6pt}}(L, X)\times \Ext_{\hspace{-2pt}\mathcal{A}\hspace{-1pt}}^1(Z, L)\to I.$ If ${\rm Hom}_{\hspace{1pt}\overline{\hspace{-2pt}\mathcal{A}\hspace{-0.6pt}}\hspace{0.6pt}}(L,X)\vspace{1pt}$ or $\Ext_{\hspace{-2pt}\mathcal{A}\hspace{-1pt}}^1(Z,L)\vspace{1pt}$ is of finite $R$-length, then it follows from Lemma \ref{duality} that $\Ext_{\hspace{-2pt}\mathcal{A}\hspace{-1pt}}^1(Z,\, L) \cong  \hspace{0.3pt} D\Hom_{\hspace{1pt}\overline{\hspace{-1pt}\mathcal{A}\hspace{-1pt}}\hspace{1pt}} (L, \,X).$ The proof of the lemma is completed.

\bigskip

The following result is a local version, but under weaker hypotheses, of the main result stated in \cite{LeZ}.

\bigskip

\begin{Theo} \label{theo_cat} Let $\mathcal{A}$ be an exact $R$-category where $R$ is artinian, and let $X, Z \in \mathcal{A}$ be strongly indecomposable with $X$ not Ext-injective and $Z$ not Ext-projective.

\vspace{-0pt}

\begin{enumerate}[$(1)$]

\item If $\,{\rm Hom}_{\hspace{1pt}\overline{\hspace{-1.5pt}\mathcal{A}\hspace{-0.6pt}}\hspace{0.6pt}}(L,X)
\in {\rm mod}\hspace{0.3pt}R\vspace{1pt}$
for all $L\in \mathcal{A}$, then $\mathcal A$ has an almost split sequence $\xymatrixcolsep{14pt}\xymatrix{0\ar[r] & X\ar[r] & Y \ar[r] & Z\ar[r]& 0}$ if and only if $\; \Ext_{\hspace{-2pt}\mathcal{A}\hspace{-1pt}}^1(Z,\, -)$ $\cong $ $\hspace{0.3pt} D\Hom_{\hspace{1pt}\overline{\hspace{-2pt}\mathcal{A}\hspace{-1pt}}\hspace{1pt}}(-, \,X).$

\vspace{2pt}

\item If $\,\Hom_{\hspace{0.3pt}\underline{\mathcal{A}\hspace{-1pt}}\hspace{1pt}}(Z,L) \in {\rm mod}\hspace{0.3pt}R\vspace{1pt}$
for all $L\in \mathcal{A}$, then $\mathcal A$ has an almost split sequence $\xymatrixcolsep{14pt}\xymatrix{0\ar[r] & X\ar[r] & Y \ar[r] & Z\ar[r]& 0}$ if and only if $\;\Ext_{\hspace{-2pt}\mathcal{A}\hspace{-1pt}}^1(-,\, X) \cong D \Hom_{\hspace{0.3pt}\underline{\mathcal{A}\hspace{-1pt}}\hspace{1pt}}(Z,\, -).$

\vspace{1pt}

\end{enumerate}\end{Theo}

\noindent{\it Proof.} We shall prove only Statement (1). Suppose that $\,{\rm Hom}_{\hspace{1pt}\overline{\hspace{-1.5pt}\mathcal{A} \hspace{-0.6pt}}}\hspace{1pt}\vspace{1pt}(L,X)$ is of finite $R$-length
for any $L\in \mathcal{A}$. Let $\phi: \Ext_{\hspace{-2pt}\mathcal{A}\hspace{-1pt}}^1(Z,\, -)\to D\Hom_{\hspace{1pt}\overline{\hspace{-2pt}\mathcal{A}\hspace{-0.6pt}}\hspace{0.6pt}}(-, X)\vspace{1pt}$ be a functorial isomorphism. Since $Z$ is not Ext-projective, $\End_{\hspace{0.5pt}\overline{\hspace{-2pt}\mathcal{A}\hspace{-0.6pt}}\hspace{0.6pt}}(X)\ne 0\vspace{1pt}$. Since $\phi_X$ is bijective, $\Ext_{\hspace{-2pt}\mathcal{A}\hspace{-1pt}}^1(Z,\, X)\vspace{1.5pt}$ is non-zero of finite $R$-length. In particular, the $\End_{\hspace{-2pt}\mathcal{A}\hspace{-1pt}}^1(X)$-socle of
$\Ext_{\hspace{-2pt}\mathcal{A}\hspace{-1pt}}^1(Z, X)$ is non-zero.
By Theorem \ref{AR-exist-cor},  $\mathcal{A}$ has a desired almost split sequence.

Conversely, let $\xymatrixcolsep{15pt}\xymatrix{0\ar[r] & X\ar[r] & Y \ar[r] & Z\ar[r]& 0}$ be an almost split sequence in $\mathcal A.$ By Theorem \ref{AR-sequence}, we have a functorial monomorphism $\phi: \Ext_{\hspace{-2pt}\mathcal{A}\hspace{-1pt}}^1(Z,\, -)\to D\Hom_{\hspace{1pt}\overline{\hspace{-2pt}\mathcal{A}\hspace{-0.6pt}}\hspace{0.6pt}}(-, X)\vspace{1pt}$. For each $L\in \mathcal{A}$, by Lemma \ref{AR-iso}(1), $\Ext_{\hspace{-2pt}\mathcal{A}\hspace{-1pt}}^1(Z,\, L) \cong \hspace{0.3pt} D\Hom_{\hspace{1pt}\overline{\hspace{-2pt}\mathcal{A}\hspace{-0.8pt}}\hspace{0.8pt}}(L, \,X),\vspace{1pt}$ and hence, the monomorphism $\phi_L: \Ext_{\hspace{-2pt}\mathcal{A}\hspace{-1pt}}^1(Z,\, L)\to D\Hom_{\hspace{1pt}\overline{\hspace{-2pt}\mathcal{A}\hspace{-0.6pt}}\hspace{0.8pt}}(L, X)\vspace{1pt}$ is an isomorphism. 
The proof of the theorem is completed.

\bigskip

We conclude this section with the following interesting consequence.

\medskip

\begin{Cor}

Let $\mathcal{A}$ be a Krull-Schmidt exact $R$-category where $R$ is artinian. If $\mathcal{A}$ has almost split sequences, then $\mathcal{A}$ is Ext-finite if and only if $\hspace{2pt}\overline{\hspace{-2pt}\mathcal{A}\hspace{-0.4pt}}\hspace{0.8pt}$ is Hom-finite, if and only if $\underline{\mathcal{A}\hspace{-1pt}}$ is Hom-finite.

\end{Cor}

\noindent {\it Proof.} Let $X, Y\vspace{1pt}\in \mathcal A$ be indecomposable. Suppose first that $\hspace{3pt}\overline{\hspace{-3pt}\mathcal A\hspace{-0.6pt}}\hspace{0.6pt}$ \hspace{1pt}is Hom-finite. If $X$ is Ext-projective, then $\Ext_{\hspace{-2pt}\mathcal{A}\hspace{-1pt}}^1(X, Y)=0\vspace{1pt}$. Otherwise, since $\Hom_{\hspace{1pt}\overline{\hspace{-2pt}\mathcal{A}\hspace{-0.6pt}}\hspace{0.6pt}}(Y, \tau_{_{\hspace{-1pt}\mathcal{A}}}X)$ is of finite $R$-length, so is $\Ext_{\hspace{-2pt}\mathcal{A}\hspace{-1pt}}^1(X, Y)\vspace{1pt}$ by Lemma \ref{AR-iso}(1). Thus $\mathcal A$ is Ext-finite.

\vspace{1.5pt}

Suppose next that $\mathcal A$ is Ext-finite. If $Y$ is Ext-injective, $\Hom_{\hspace{1pt}\overline{\hspace{-1.5pt}\mathcal{A}\hspace{-0.6pt}}\hspace{0.6pt}}(X, Y)= 0\vspace{1.5pt}$. Otherwise, since $\Ext_{\hspace{-2pt}\mathcal{A}\hspace{-1pt}}^1(\tau_{_{\hspace{-1pt}\mathcal{A}}}^-Y, X)\vspace{2pt}$ is of finite $R$-length, so is   $\Hom_{\hspace{1pt}\overline{\hspace{-2pt}\mathcal{A}\hspace{-0.6pt}}\hspace{0.6pt}}(X,Y)$ by Lemma \ref{AR-iso}(1). This shows that $\hspace{3pt}\overline{\hspace{-2pt}\mathcal{A}\hspace{-0.4pt}}\hspace{0.8pt}$ is Hom-finite. Similarly, $\mathcal A$ is Ext-finite if and only if $\hspace{1pt}\underline{\mathcal A\hspace{-1.5pt}}\hspace{1.5pt}$ is Hom-finite. The proof of the corollary is completed.


\smallskip

{\center \section{Minimal approximations}}

\smallskip

Throughout this section, $\mathcal{A}$ stands for an exact $R$-category, and $\mathcal{C}$ for an exact subcategory of $\mathcal{A}$, that is, $\mathcal C$ is an extension-closed full subcategory of $\mathcal{A}$. The objective of this section is to study when an almost split sequence in $\mathcal A$ induces an almost split sequence in $\mathcal C$. 

\bigskip

We start with some notation and terminology. Let $\widetilde{\mathcal{C}\hspace{-0.5pt}}$ and $\underset{\widetilde{}}{\mathcal{C}}\vspace{-4pt}$ stand for the full subcategories generated by the objects in $\mathcal C$ of $\hspace{2pt}\overline{\hspace{-3pt}\mathcal A\hspace{-0.6pt}}\hspace{0.6pt}$  and  $\underline{\mathcal A \hspace{-1pt}}\hspace{1.5pt}$, respectively.
Fix an object $X\in \mathcal{A}$. A morphism $f: M\to X$ in $\mathcal A$ with $M\in \mathcal{C}$ is called a {\it right injectively stable $\mathcal C$-approximation} of $X$ if $\bar f$ is a right $\widetilde{\mathcal{C}\hspace{-0.5pt}}$-approximation of $X$ in $\hspace{3pt}\overline{\hspace{-3pt} \mathcal A \hspace{-0.5pt}}\,$; and a {\it minimal right injectively stable $\mathcal C$-approximation} if, in addition, $\bar f$ is right minimal in $\hspace{3pt}\overline{\hspace{-3pt} \mathcal A \hspace{-0.5pt}}\,$ and $M$ has no non-zero summand which is Ext-injective in $\mathcal A$. Dually, a morphism $g: X\to N$ in $\mathcal A$ with $N\in \mathcal{C}$ is called a {\it left projectively stable $\mathcal C$-approximation} of $X$ if $\underline{\hspace{1pt}g\hspace{-2pt}}\hspace{2pt}$ is a left $\underset{\widetilde{}}{\mathcal{C}}\vspace{-2pt}$-approximation of $X$ in $\underline{\hspace{1pt} \mathcal{A}\hspace{-0.5pt}}\,$; and a {\it minimal left projectively stable $\mathcal C$-approximation} if, in addition, $\underline{\hspace{1pt}g\hspace{-1.5pt}}\hspace{2pt}\vspace{1pt}$ is left minimal in $\underline{\hspace{1pt}\mathcal A \hspace{-0.5pt}}\,$ and $N$ has no non-zero summand which is Ext-projective in $\mathcal A$.

\bigskip

\begin{Lemma}\label{approx-mini}

Let $X\in \mathcal A$ have a right injectively stable $\mathcal C$-aprroximation $f: M\to X.$ If $M$ is Krull-Schmidt, 
then $f$ decomposes as $f=(g, h): N\oplus L\to X$, where $g$ is a minimal right injectively stable $\mathcal C$-approximation of $X$.

\end{Lemma}

\noindent{\it Proof.} Let $M$ be Krull-Schmidt. Then $f=(f_1, \ldots, f_r): M=M_1\oplus \cdots \oplus M_r\to X,$ where the $M_i$ are strongly indecomposable in $\mathcal C$. If $f$ is injectively trivial in $\mathcal A$, then $0:0\to X$ is a minimal injectively stable $\mathcal C$-approximation of $X$. Otherwise, we may assume that there exists some $1\le s\le r$ such that $M_i$ is Ext-injective in $\mathcal A$ if and only if $s<i\le r.$ Put $U=M_1\oplus \cdots \oplus M_s$ and $\vspace{0pt}u=(f_1, \ldots, f_s): U\to X$. Then $\bar{u\hspace{0.3pt}}$ is a right $\widetilde{\mathcal{C}\hspace{-0.3pt}}$-approximation of $X$ in $\hspace{2pt}\overline{\hspace{-3pt}\mathcal A\hspace{-0.6pt}}\hspace{0.6pt}$. By Lemma \ref{KS}, the $M_i$ with $1\le i\le s\vspace{1pt}\,$ are strongly indecomposable in $\hspace{2pt}\overline{\mathcal{\hspace{-3pt}A\hspace{-0.6pt}}}\hspace{0.6pt}$. Hence, by Proposition \ref{KS-idemp},
the idempotents in $\End_{\hspace{1pt}\overline{\mathcal{\hspace{-2pt}A\hspace{-1pt}}}\hspace{0.6pt}}(U)\vspace{1.5pt}$ split in $\hspace{2pt}\overline{\mathcal{\hspace{-3pt}A\hspace{-0.6pt}}}\hspace{0.6pt}$. Therefore,
there exists a decomposition
$\overline{\hspace{-0.3pt}u\hspace{-0.5pt}}=(\bar{v}, \bar{0}): U=V\oplus W \to X,$
where $\bar{v}$ is right minimal in $\hspace{1pt}\overline{\mathcal{\hspace{-2pt}A\hspace{-0.6pt}}}\hspace{0.6pt};$
see \cite[(1.4)]{KSa}. Then $\bar{v}: V\to X$ is a minimal right $\widetilde{\mathcal{C}\hspace{-0.3pt}}$-approximation of $X$ in $\hspace{2pt}\overline{\mathcal{\hspace{-3pt}A\hspace{-0.6pt}}}\hspace{0.6pt}$. Since $V$ is a direct summand of $M$ in $\hspace{2pt}\overline{\mathcal{\hspace{-3pt}A\hspace{-0.6pt}}}\hspace{0.6pt}$, by Proposition \ref{KS-idemp}, we may assume that $\hspace{2pt}\overline{\mathcal{\hspace{-3pt}A\hspace{-0.6pt}}}\hspace{0.6pt}$ has an isomorphism $\bar{p}: N\to V$, where $N=M_1\oplus \cdots \oplus M_t$ for some $1\le t\le s$. Setting $g=vp$, we get a minimal right
$\widetilde{\mathcal{C}\hspace{-0.3pt}}$-approximation $\bar{g}: N\to X$ of
$X$ in $\hspace{2pt}\overline{\mathcal{\hspace{-3pt}A\hspace{-0.6pt}}}\hspace{0.6pt}$. Then, $g: N\to X$ is a minimal right injectively stable $\mathcal{C}$-approximation of $X$. The proof of the lemma is completed.

\bigskip

The following result characterizes the Ext-projective or Ext-injective objects in $\mathcal C$ which admit an almost split sequence in $\mathcal A$.

\medskip

\begin{Lemma}\label{Ext-proj}

Let $\xymatrixcolsep{14pt}\xymatrix{0\ar[r] & X\ar[r] & Y \ar[r] & Z\ar[r]& 0}$ be an almost split sequence in $\mathcal A$.

\vspace{-1pt}

\begin{enumerate}[$(1)$]

\item If $Z\in {\mathcal C}\vspace{1pt}$, then it is Ext-projective in $\mathcal C$ if and only if $\;0\to X$ is a right injectively stable $\mathcal C$-approximation of $X$.

\vspace{1pt}

\item If $X\in {\mathcal C}$, then it is Ext-injective in $\mathcal C$ if and only if $Z\to 0$ is a left projectively stable $\mathcal C$-approximation of $Z$.


\end{enumerate} \end{Lemma}

\noindent{\it Proof.} We shall prove only Statement (1). For each $L\in {\mathcal C}\vspace{1pt}$, by Proposition \ref{bilinearform}, there exists a non-degenerate $R$-bilinear form $<\; , \;>: \Hom_{\hspace{1pt}\overline{\hspace{-2pt}\mathcal{A}\hspace{-1pt}}\hspace{1pt}}(L, X)\times \Ext_{\hspace{-2pt}\mathcal{A}\hspace{-1pt}}^1(Z, L) \vspace{1pt}\to I$. If $Z\in {\mathcal C}\vspace{1pt}$, then $\Ext_{\hspace{-2pt}\mathcal{A}\hspace{-1pt}}^1(Z, L)=\Ext^1_{\mathcal C}(Z, L)$. Therefore, $Z$ is Ext-projective in $\mathcal C$ if and only if $\Ext^1_{\mathcal C}(Z, L)=0$ for all $L\in {\mathcal C}$, if and only if $\Hom_{\hspace{1pt}\overline{\hspace{-1.5pt}\mathcal{A}\hspace{-1pt}}\hspace{1pt}}(L, X)=0\vspace{1pt}$ for all $L\in \mathcal{C},$ that is, $\;0\to X$ is a right injectively stable $\mathcal C$-approximation of $X$. The proof of the lemma is completed.

\bigskip

We shall now show that minimal injectively or projectively
stable $\mathcal{C}$-approxima\-tions of almost split sequences in $\mathcal{A}$ are almost split sequences in $\mathcal{C}$. For this purpose,  we need the following preparatory lemma.

\medskip

\begin{Lemma}\label{conseq-approx}

Let $\,\delta: \xymatrixcolsep{14pt}\xymatrix{0\ar[r] & X\ar[r] & Y \ar[r] & Z\ar[r]& 0}$ be an almost split sequence in $\mathcal A$, where $Z$ lies in $\mathcal{C}$ and $X$ admits a minimal right injectively stable $\mathcal C$-approximation $f: M\to X$.

\vspace{-1pt}

\begin{enumerate}[$(1)$]

\item For $L\in {\mathcal C}$, the map $\Ext_{\hspace{-2pt}\mathcal{A}\hspace{-1pt}}^1(L, f): \Ext_{\hspace{-2pt}\mathcal{A}\hspace{-1pt}}^1(L, M) \to \Ext_{\hspace{-2pt}\mathcal{A}\hspace{-1pt}}^1(L, X)$ is injective.

\vspace{2pt}

\item If $u: L\to M$ lies in $\mathcal C$, then $fu\in I_{\hspace{-1pt}\mathcal A}(L, X)$ if and only if $u\in I_{\mathcal C\hspace{0.2pt}}(L, M)$.

\vspace{2pt}

\item If $Z$ is not Ext-projective in $\mathcal C$, then $M$ is indecomposable.

\end{enumerate}
\end{Lemma}

\noindent {\it Proof.} (1) Let $L\in \mathcal{C}$, and consider a pushout diagram \vspace{-1pt}
$$\xymatrixcolsep{16pt}\xymatrixrowsep{16pt}\xymatrix{
\eta: \;   & 0\ar[r] & M\ar[r]^r \ar[d]_f & E \ar[r] \ar[d]^g & L\ar[r]\ar@{=}[d]& 0\\
\zeta: \;  & 0\ar[r] & X\ar[r]^s & F \ar[r] & L\ar[r]& 0} \vspace{-2pt} $$ in $\mathcal A$. If
$\zeta$ splits, then $ts=\id$ for some $t: F\to X$, and hence $f=(tg)r$. Since $E\in \mathcal{C}$, there exists some $h: E\to M$ such that $\overline{tg}=\overline{fh}.$ This gives rise to  $\bar f=\bar f \cdot \overline{hr}.\,$ Since $\bar f$ is right minimal in $\hspace{3pt}\overline{\hspace{-3pt} \mathcal A \hspace{-0.5pt}}\,$, we see that $\overline{hr}$ is an automorphism of $M$ in $\hspace{3pt}\overline{\hspace{-3pt} \mathcal A \hspace{-0.5pt}}\,$. Therefore, $\overline{wr}=\bar \id$, for some $w: E\to M$. That is, $\id-wr$ is injectively trivial in $\mathcal A$. In particular, $\id-wr$ factors through $r$. Then $r$ is a section, that is, $\eta$ splits.

\vspace{1pt}

(2) Let $u: L\to M$ be a morphism in $\mathcal C$. Assume first that $fu$ is not injectively trivial in $\mathcal A$. Let $\varphi: \Ext_{\hspace{-2pt}\mathcal{A}\hspace{-1pt}}^1(Z, X)\to I$ be an $R$-linear form such that $\varphi(\delta)\ne 0$. In view of Proposition \ref{bilinearform}, there exists $\zeta\in \Ext_{\hspace{-2pt}\mathcal{A}\hspace{-1pt}}^1(Z, L)$ such that $\varphi(fu\hspace{0.4pt}\zeta)\ne 0$. In particular, $u\,\zeta\ne 0$. Since $\zeta$ lies in $\mathcal C$, we see that $u\not\in I_{\mathcal C}(L, M)$. Suppose conversely that $fu$ is injectively trivial in $\mathcal A$. If $\eta: \xymatrixcolsep{15pt}\xymatrix{0\ar[r] & L\ar[r] & E \ar[r] & N\ar[r]& 0}$ is a short exact sequence in $\mathcal C$, then $(fu)\eta=0$, that is, $f(u\eta)=0$. By Statement (1), $u\eta=0$. Therefore, $u\in I_{\mathcal C}(L, M)$.

\vspace{2pt}

(3) Since the right $\End_{\hspace{-1pt}\mathcal{A}\hspace{-1pt}}(Z)\vspace{1pt}$-module $\Ext_{\hspace{-2pt}\mathcal{A}\hspace{-1pt}}^1(Z, X)$ has a simple socle, every non-zero $\End_{\hspace{-1pt}\mathcal{A}\hspace{-1pt}}(Z)$-submodule of $\Ext_{\hspace{-2pt}\mathcal{A}\hspace{-1pt}}^1(Z, X)\vspace{1pt}$ is indecomposable.
Suppose that $Z$ is not Ext-projective in $\mathcal C$. In particular, it is not Ext-projective in $\mathcal A$. By Lemma \ref{Ext-proj}(1), $f$ is not injectively trivial in $\mathcal A$. By Proposition \ref{bilinearform},
$\Ext_{\hspace{-2pt}\mathcal{A}\hspace{-1pt}}^1(Z, M) \vspace{1pt}\ne 0$. Observe that $\Ext_{\hspace{-2pt}\mathcal{A}\hspace{-1pt}}^1(Z, f): \Ext_{\hspace{-2pt}\mathcal{A}\hspace{-1pt}}^1(Z, M) \to \Ext_{\hspace{-2pt}\mathcal{A}\hspace{-1pt}}^1(Z, X)\vspace{1.5pt}$ is $\End_{\hspace{-1pt}\mathcal{A}\hspace{-1pt}}(Z)$-linear and injective by Statement (1). Thus,
$\Ext_{\hspace{-2pt}\mathcal{A}\hspace{-1pt}}^1(Z, M) \vspace{1pt}$ is isomorphic to a non-zero $\End_{\hspace{-1pt}\mathcal{A}\hspace{-1pt}}(Z)$-submodule of $\Ext_{\hspace{-2pt}\mathcal{A}\hspace{-1pt}}^1(Z, X)\vspace{1pt}$, and hence, it is an indecomposable right
$\End_{\hspace{-1pt}\mathcal{A}\hspace{-1pt}}(Z)$-module. Assume that $M=M_1\oplus M_2,\vspace{1pt}$ with non-zero injections $q_i: M_i\to M$, $i=1, 2.$ By the hypothesis, \vspace{1pt} $M_1$ and $M_2$
are non-zero in $\hspace{2pt}\overline{\hspace{-2pt}\mathcal{A}\hspace{-0.4pt}}\hspace{0.4pt}$. Since $\hspace{1pt}\overline{\hspace{-1pt}f}$ is right minimal in $\hspace{2pt}\overline{\hspace{-2pt}\mathcal A\hspace{-0.6pt}}\hspace{0.6pt}$, we have
$\overline{fq_i\hspace{-2pt}}\ne \bar{0}\vspace{2pt}$, and by Proposition \ref{bilinearform}, $\Ext_{\hspace{-2pt}\mathcal{A}\hspace{-1pt}}^1(Z, M_i)\ne 0\vspace{1pt}$, $i=1, 2.$ This is absurd, since $\Ext_{\hspace{-2pt}\mathcal{A}\hspace{-1pt}}^1(Z, M)\cong \Ext^1_{\mathcal{A}}(Z, M_1) \oplus \Ext^1_{\mathcal{A}}(Z, M_2)\vspace{1pt}$. The proof of the lemma is completed.

\bigskip



\begin{Prop}\label{approx-ass}

Let $\,\delta: \xymatrixcolsep{14pt}\xymatrix{0\ar[r] & X\ar[r] & Y \ar[r] & Z\ar[r]& 0}$ be an almost split sequence in $\mathcal A$, where $Z$ lies in ${\mathcal C}$ and $X$ has a non-zero minimal right injectively stable $\mathcal C$-approxi\-mation $f: M\to X$. If $\,M$ is Krull-Schmidt, then $\mathcal{A}$ has a pushout diagram
\vspace{-6pt} $$\xymatrixcolsep{14pt}\xymatrixrowsep{16pt}\xymatrix{
\eta: & 0\ar[r] & M\ar[r]\ar[d]_f & N \ar[d] \ar[r] & Z\ar[r]\ar@{=}[d]& 0\\
\delta: & 0\ar[r] & X\ar[r] & Y \ar[r] & Z\ar[r]& 0;
}\vspace{1pt}$$ and in any such pushout diagram, $\eta$ is an almost split sequence in $\mathcal C$.

\end{Prop}

\noindent{\it Proof.} By assumption, $M$ is non-zero. Then, by definition, $M$ is not Ext-injective in $\mathcal A$. That is, $M$ is non-zero in $\hspace{2pt}\overline{\hspace{-3pt}\mathcal A\hspace{-0.6pt}}\hspace{0.6pt}$. Being right minimal, $\bar{f}$ is non-zero in $\hspace{2pt}\overline{\hspace{-3pt}\mathcal A\hspace{-0.6pt}}\hspace{0.6pt}$. By Lemma \ref{Ext-proj}(1), $Z$ is not
Ext-projective in $\mathcal{C}$. Then, $M$ is indecomposable by Lemma \ref{conseq-approx}(3), and hence strongly indecomposable since it is Krull-Schmidt. For each $L\in {\mathcal C},$ in view of Lemma \ref{conseq-approx}(2), we have an isomorphism
$$\Hom_{\hspace{0.5pt}\overline{\mathcal{C}\hspace{-0.5pt}}}\hspace{0.8pt}(L, M) \to \Hom_{\hspace{1pt}\overline{\hspace{-2pt}\mathcal{A}\hspace{-0.6pt}}\hspace{0.6pt}}(L, X): \tilde{u}\mapsto \overline{\hspace{-1pt}fu\hspace{-1pt}}\hspace{1pt},$$ which is clearly natural in $L$. This induces a functorial isomorphism
$$D\Hom_{\hspace{1pt}\overline{\hspace{-2pt} \mathcal{A}\hspace{-0.6pt}}\hspace{0.6pt}}(-, X) |_{\mathcal C} \to D\Hom_{\hspace{0.5pt}\overline{\mathcal{C} \hspace{-0.5pt}}}\hspace{0.5pt}(-, M).$$
On the other hand, by Theorem \ref{AR-sequence}, there exists a functorial monomorphism
$$\Ext_{\hspace{-2pt}\mathcal{A}\hspace{-1pt}}^1(Z, -)\to D\Hom_{\hspace{1pt}\overline{\hspace{-2pt}\mathcal A\hspace{-0.6pt}}\hspace{0.6pt}}(-, X).$$
Since $\Ext_{\mathcal C}^1(Z, -)=\Ext_{\hspace{-2pt}\mathcal{A}\hspace{-1pt}}^1(Z, -)|_{\mathcal C},\vspace{3pt}$ we get a functorial monomorphism
$$\phi: \Ext_{\mathcal C}^1(Z, -)\to D\Hom_{\hspace{0.5pt}\overline{\mathcal{C}\hspace{-0.6pt}}\hspace{0.6pt}}(-, M).$$

Since $f$ is not injectively trivial, by Proposition \ref{bilinearform}, there exists a non-zero extension $\zeta\in \Ext_{\hspace{-2pt}\mathcal{A}\hspace{-1pt}}^1(Z, M)$ such that $f\zeta\ne 0$.
Since $\delta$ is almost split in $\mathcal A$, there exists $g\in \End_{\hspace{-0.5pt}\mathcal{A}\hspace{-1pt}}(Z)$ such that $\delta=(f\zeta)g=f(\zeta g)$. This establishes the existence of a commutative diagram as stated in the proposition.

\vspace{1pt}

Now, let $\eta\in \Ext^1_{\mathcal C}(M, Z)$ be such that $f\eta=\delta$. Suppose that $u\eta\ne 0$ for some $u\in {\rm rad}(\End_{\hspace{0.5pt}\mathcal{C}\hspace{-1pt}}(M)).$ Since $\delta$ is almost split, there exists some $v: M\to X$ such that $v(u\eta)=\delta$.
Since $\bar{f}$ is a right $\widetilde{\mathcal{C}\hspace{-0.5pt}}$-approximation of $X$, there exists some $w: M\to M$ such that $\overline{v}=\overline{fw}$. This yields $f\eta=\delta=f(wu\eta).$ By Lemma \ref{conseq-approx}(1), $\eta=(wu)\eta$.
Since $wu\in  {\rm rad}(\End_{\hspace{0.5pt}\mathcal{C}\hspace{-1pt}}(M))$, we get $\eta=0$, a contradiction. This proves that
$\eta$ lies in the $\End_{\mathcal{C}\hspace{-0.5pt}}(M)$-socle of $\Ext_{\hspace{-2pt}\mathcal{A}\hspace{-1pt}}^1(M, Z)$. Since $\phi_M$ is an $\End_{\mathcal{C}\hspace{-0.5pt}}(M)$-linear monomorphism, $\phi_M(\eta)$ is almost vanishing on $\End_{\hspace{0.5pt}\overline{\mathcal{C}\hspace{-0.5pt}}}(M)$. By Theorem \ref{AR-sequence}, $\eta$ is an almost split sequence in $\mathcal{C}$. The proof of the proposition is completed.

\bigskip

The following statement generalizes the main results stated in \cite{AuS}; see also \cite{GaR}. Observe that we do not impose any finiteness assumption.

\medskip

\begin{Cor}\label{ass-subcat}

Let $\mathcal A$ have almost split sequences. If $\mathcal C$ is Krull-Schmidt and functorially finite in $\mathcal A$, then $\mathcal C$ has almost split sequences.

\end{Cor}

\noindent{\it Proof.} Assume that $\mathcal C$ is Krull-Schmidt and functorially finite in $\mathcal A$. Let $Z\in\mathcal{C}$ be indecomposable but not Ext-projective. Then $\mathcal A$ has an almost split sequence $\xymatrixcolsep{15pt}\xymatrix{0\ar[r] & X\ar[r] & Y \ar[r] & Z\ar[r]& 0.}$ By the assumption and Lemma \ref{approx-mini}, $X$ has a minimal right injectively stable $\mathcal C$-approximation $f: M\to X$. By Lemma \ref{Ext-proj}(1), $f$ is non-zero. By Proposition \ref{approx-ass}, $\mathcal C$ has a almost split sequence ending with $Z$. This shows that $\mathcal C$ has right almost split sequences. Dually, $\mathcal C$ has left almost split sequences. The proof of the corollary is completed.

\bigskip

Next, we shall establish the converse of Proposition \ref{approx-ass}. For this purpose, some finiteness assumption is needed.

\medskip

\begin{Prop}\label{ass-approx}

Let $R$ be artinian, and let $\,\delta: \xymatrixcolsep{14pt}\xymatrix{0\ar[r] & X\ar[r] & Y \ar[r] & Z\ar[r]& 0}$ and $\,\eta: \xymatrixcolsep{15pt}\xymatrix{0\ar[r] & M\ar[r] & N \ar[r] & Z\ar[r]& 0}\vspace{1pt}$ be almost split sequences in $\mathcal A$ and $\mathcal C$, respectively. If $\,\Hom_{\hspace{0.2pt}\overline{\mathcal{C}\hspace{-0.5pt}}}\hspace{0.5pt}(L, M)
\in {\rm mod}\hspace{0.3pt}R$
for any $L\in \mathcal{C},$
then $\eta$ embeds in a pushout diagram
$$\xymatrixcolsep{14pt}\xymatrixrowsep{16pt}\xymatrix{
\eta: & 0\ar[r] & M\ar[r]\ar[d]_f & N \ar[d]^g \ar[r] & Z\ar[r]\ar@{=}[d] & 0\\
\delta: & 0\ar[r] & X\ar[r] & Y \ar[r] & Z\ar[r]& 0
}$$ in $\mathcal{A};$ and in any such pushout diagram, $f$ is a minimal right injectively stable $\mathcal C$-approxima\-tion of $X,$ and $g$ is a right injectively stable ${\mathcal C}$-approximation of $\hspace{0.5pt}Y.$

\end{Prop}

\noindent{\it Proof.} Suppose that $\Hom_{\hspace{0.5pt}\overline{\hspace{-0.3pt}\mathcal{C}\hspace{-0.5pt}}} \hspace{1pt}(L, M)$ is of finite $R$-length, for every $L\in \mathcal{C}$. Since $\delta$ is almost split, $\mathcal A$ has a commutative diagram
$$\xymatrixcolsep{16pt}\xymatrixrowsep{16pt}\xymatrix{
\eta: & 0\ar[r] & M\ar[r]^r\ar[d]_f & N \ar[d]^g \ar[r]^s & Z\ar[r]\ar@{=}[d]& 0\\
\delta: & 0\ar[r] & X\ar[r]^u & Y \ar[r]^v & Z\ar[r]& 0.
} \vspace{1pt}$$

Fix an $R$-linear form $\varphi: \Ext_{\hspace{-2pt}\mathcal{A}\hspace{-1pt}}^1(Z, X)\to I$ such that $\varphi(\delta)\ne 0$. This yields an $R$-linear form
\vspace{-1pt}$$\psi: \Ext^1_{\mathcal C}(Z, M) \to I: \zeta \mapsto \varphi(f \zeta)
\vspace{1pt}$$ such that $\psi(\eta)=\varphi(\delta)\ne 0$.
Let $L\in \mathcal{C}$. By Proposition \ref{bilinearform}, we have non-degenerate $R$-bilinear forms
$$<\, , \, >_\varphi: \; \Hom_{\hspace{0.5pt}\overline{\hspace{-2pt}\mathcal{A}\hspace{-0.6pt}}\hspace{0.5pt}}(L, X)\times \Ext_{\hspace{-2pt}\mathcal{A}\hspace{-1pt}}^1(Z, L)\to I: (\overline{\hspace{-0.5pt}q\hspace{-1.5pt}}\,, \zeta) \mapsto \varphi(q\zeta) \vspace{-5pt}$$
and
$$<\, , \, >_\psi: \; \Hom_{\overline{\mathcal{C}\hspace{-0.8pt}}\hspace{0.8pt}} (L, M)\times \Ext^1_{\mathcal C}(Z, L)\to I: (\tilde{p}\,, \zeta) \mapsto \psi(p\zeta). \vspace{2pt}$$

Let $q: L\to X$ be a morphism in $\mathcal A$. Since $\Ext^1_{\hspace{-1pt}\mathcal{A}}(Z, L)=\Ext^1_{\mathcal C}(Z, L)$, we have an $R$-linear form \vspace{-2pt}
$$<\hspace{-2pt}\bar q, \, -\hspace{-3pt}>_\varphi\hspace{0.5pt}: \; \Ext^1_{\mathcal C}(Z, L)\to I: \zeta\mapsto <\hspace{-2pt}\bar q, \, \zeta\hspace{-3pt}>_\varphi.
\vspace{2pt}$$

Since $\Hom_{\hspace{0.5pt}\overline{\hspace{-0pt}\mathcal{C}\hspace{-1pt}}} \hspace{1pt}(L, M)$ is of finite $R$-length, by Lemma \ref{duality}, there exists $p: L\to M$ in $\mathcal C$ such that $<\hspace{-2pt}\bar q, \; -\hspace{-3pt}>_\varphi  \, =\, <\hspace{-2pt}\tilde p\,, \; - \hspace{-3pt}>_\psi$. That is, for any $\zeta\in \Ext_{\hspace{-2pt}\mathcal{A}\hspace{-1pt}}^1(Z, L)$, we have
$$<\hspace{-3pt}\bar q, \, \zeta\hspace{-2.5pt}>_\varphi\,=\,<\hspace{-3pt}\tilde p\,, \, \zeta\hspace{-2pt}>_\psi=\psi(p\,\zeta)=\varphi(fp\,\zeta)=<\hspace{-2pt}\overline{\hspace{-1pt}fp\hspace{-0.5pt}}\,, \zeta\hspace{-2pt}>_\varphi.$$

Since $<\, , \, >_\varphi$ is non-degenerate, $\bar q=\overline{\hspace{-1pt}fp\hspace{0.3pt}}\,$. This shows that $f$ is a right injectively stable $\mathcal C$-approximation of $X$, which is minimal since $M$ is strongly indecomposable in both $\mathcal A$ and $\hspace{2pt}\overline{\hspace{-2pt}\mathcal A\hspace{-0.6pt}}\hspace{0.6pt}$; see (\ref{KS}).

Next, consider a morphism $h: L\to Y$ in $\mathcal A$ with $L \in \mathcal{C}$. Since $vh$ is not a retraction, there exists $w: L\to N$ such that $vh=sw=vgw$. Thus $h-gw=ut$ for some $t: L\to X$. Using what we have just proved, there exists some morphism $j: L\to M\vspace{1pt}$ such that $\bar{t}=\overline{fj}$, and hence $\overline{h-gw}=\overline{ufj}=\overline{grj}.$ This yields that
$\overline{h}=\bar{g} (\overline{w}+\overline{rj})$. That is, $g$ is a right injectively stable $\mathcal C$-approximation of $Y$.
The proof of the proposition is completed.

\bigskip

We are ready to obtain the main result of this section.

\medskip

\begin{Theo}\label{subcat-exist}

Let $\mathcal{A}$ be an exact $R$-category where $R$ is artinian, and let $\mathcal{C}$ be an exact subcategory of $\mathcal A$ which is Ext-finite and Krull-Schmidt.

\vspace{-1pt}

\begin{enumerate}[$(1)$]

\item If $\mathcal A$ has right almost split sequences, then $\mathcal{C}$ has right almost split sequences if and only if $\,\tau_{_{\hspace{-2pt}\mathcal A}}\hspace{-1pt}Z$ has a right injectively stable $\mathcal C$-approximation, for any indecomposable not Ext-projective object $Z$ in $\mathcal C$.

\vspace{1pt}

\item If $\mathcal A$ has left almost split sequences, then $\mathcal{C}$ has left almost split sequences if and only if $\,\tau_{_{\hspace{-2pt}\mathcal A}}\hspace{-1pt}^{\hspace{-2pt}-\hspace{-1pt}}X\vspace{1pt}$ has a left projectively stable $\mathcal C$-approximation, for any indecomposable not Ext-injective object $X$ in $\mathcal C$.

\end{enumerate}
\end{Theo}

\noindent{\it Proof.} We shall prove only Statement (1). Assume that $\mathcal A$ has right almost split sequences.
Since $\mathcal C$ is Krull-Schmidt, as seen in the proof of Corollary \ref{ass-subcat}, the sufficiency follows from Proposition \ref{approx-ass} and Lemmas \ref{approx-mini} and \ref{Ext-proj}. For proving the necessity, let $Z\in \mathcal{C}$ be indecomposable such that $\mathcal C$ has an almost split sequence $\xymatrixcolsep{15pt}\xymatrix{0\ar[r] & M\ar[r] & N \ar[r] & Z\ar[r]& 0.}$ For any $L\in \mathcal{C}$, since
$\Ext^1_{\mathcal C}(Z, L)$ is of finite $R$-length, so is  $\,\Hom_{\hspace{0.2pt}\overline{\mathcal{C}\hspace{-0.5pt}}}\hspace{0.5pt}(L, M)$
by Lemma \ref{AR-iso}(1). Thus, by Proposition \ref{ass-approx}, $\,\tau_{_{\hspace{-2pt}\mathcal A}}\hspace{-1pt}Z$ has a minimal right injectively stable $\mathcal C$-approximation. The proof of the theorem is completed.

\bigskip

We conclude this section with an application. Following \cite{Dic}, we shall say that a pair $({\mathcal T}, {\mathcal F})$ of full subcategories of $\mathcal A$ is a {\it torsion theory} if the following two conditions are satisfied$\,:$

\vspace{-1pt}

\begin{enumerate}[$(1)$]

\item For any objects $T\in {\mathcal T}$ and $F\in {\mathcal F}$, we have $\Hom_{\hspace{-1pt}\mathcal{A}\hspace{-1pt}}(T, F)=0.$

\vspace{1pt}

\item For any object $X\in \mathcal A$, there exists a canonical short exact sequence
\vspace{-3pt}
$$\xymatrixcolsep{14pt}\xymatrix{0\ar[r] & t(X)\ar[r] & X \ar[r] & f(X)\ar[r]& 0,}$$
where $t(X)\in {\mathcal T}$ and $f(X)\in {\mathcal F}$.

\end{enumerate}

\bigskip

{\sc Remark.} If $({\mathcal T}, {\mathcal F})$ is a torsion theory in $\mathcal A$,  then it is easy to see that ${\mathcal T}$ and ${\mathcal F}$ are exact subcategories of $\mathcal A$.

\bigskip



\begin{Prop}

Let  $\mathcal A$ be an exact $R$-category with a torsion theory $({\mathcal T}, {\mathcal F}),$ and let  $\xymatrixcolsep{14pt}\xymatrix{0\ar[r] & X\ar[r] & Y \ar[r] & Z\ar[r]& 0}$ be an almost split sequence in $\mathcal A$.

\begin{enumerate}[$(1)$]

\item If $Z\in \mathcal{T}$ is not Ext-projective, then $\mathcal T$ has an induced almost split sequence $\xymatrixcolsep{14pt}\xymatrix{0\ar[r] & t(X)\ar[r] & t(Y) \ar[r] & Z\ar[r]& 0.}$

\item If $X\in {\mathcal F}$ is not Ext-injective, then $\mathcal F$ has an induced almost split sequence $\xymatrixcolsep{14pt}\xymatrix{0\ar[r] & X\ar[r] & f(Y) \ar[r] & f(Z)\ar[r]& 0.}$

\end{enumerate}

\end{Prop}

\noindent{\it Proof.} We shall only prove Statement (1). Suppose that $Z\in \mathcal T$ is not Ext-projective. Consider the canonical short exact sequence \vspace{-1pt}
$$\xymatrixcolsep{16pt}\xymatrix{0\ar[r] & t(X)\ar[r]^-q & X \ar[r]^-p & f(X)\ar[r]& 0.}\vspace{-2pt}$$
Observe that $q$ is a right $\mathcal T$-approximation, and hence a right injectively stable $\mathcal T$-approximation, of $X$. By Lemma \ref{Ext-proj}(1), $q$ is not injectively trivial in $\mathcal A$. In particular, $t(X)$ is not Ext-injective in $\mathcal A$. Moreover, since $q$ is a mono\-morphism, we have an $R$-linear isomorphism
$$\Hom_{\mathcal{A}}(t(X), q): \End_{\hspace{-1pt}\mathcal{A}\hspace{-1pt}}(t(X))\to \Hom_{\mathcal A}(t(X), X).$$
On the other hand, $\Hom_{\hspace{0.5pt}\underline{\mathcal{A}\hspace{-1pt}}\hspace{1pt}}(Z, f(X))=0\vspace{1pt}$ since $Z\in {\mathcal T}$, and thus,
$\Ext_{\hspace{-2pt}\mathcal{A}\hspace{-1pt}}^1(f(X), X)=0$ by Proposition \ref{bilinearform}. Applying $\Hom_{\mathcal A}(-, X)$ to the canonical short exact sequence yields an $R$-linear epimorphism
$$\Hom_{\mathcal A}(q, X): \End_{\hspace{-1pt}\mathcal{A}\hspace{-1pt}}(X)\to \Hom_{\mathcal A}(t(X), X).$$ Composing this with the inverse of $\Hom_{\mathcal A}(t(X), q)$, we get an $R$-linear epimorphism
$\varphi: \End_{\hspace{-1pt}\mathcal{A}\hspace{-1pt}}(X)\to \End_{\hspace{-1pt}\mathcal{A}\hspace{-1pt}}(t(X)),$ which is a ring morphism. Since $\End_{\hspace{-1pt}\mathcal{A}\hspace{-1pt}}(X)$ is local, so is $\End_{\hspace{-1pt}\mathcal{A}\hspace{-1pt}}(t(X))$. Hence, $q: t(X)\to X$ is a minimal right injectively stable $\mathcal T$-approximation of $X$. By Proposition \ref{approx-ass}, $\mathcal A$ has a pushout diagram
\vspace{-0pt} $$\xymatrixcolsep{16pt}\xymatrixrowsep{16pt}\xymatrix{
0\ar[r] & t(X)\ar[r]\ar[d]_q & N \ar[d]^g \ar[r] & Z\ar[r]\ar@{=}[d]& 0\\
0\ar[r] & X\ar[r] & Y \ar[r] & Z\ar[r]& 0,
}\vspace{-0pt}$$ where the upper row is an almost split sequence in $\mathcal T$. Using the Snake Lemma, we infer that
$N\cong t(Y)$. The proof of the proposition is completed.

\bigskip

The following statement is an immediate consequence of the preceding result.

\medskip

\begin{Cor}

Let  $\mathcal A$ be an exact $R$-category with a torsion theory $({\mathcal T}, {\mathcal F})$.

\vspace{-2pt}

\begin{enumerate}[$(1)$]

\item If $\mathcal A$ has right almost split sequences, then $\mathcal T$
has right almost split sequences.

\item If $\mathcal A$ has left almost split sequences, then $\mathcal F$  has left almost split sequences.

\end{enumerate}
\end{Cor}

\smallskip

{\center \section{Representations of infinite quivers}}

\smallskip

Let $k$ be a field, and let $Q=(Q_0, Q_1)$ be a quiver which is locally finite and interval-finite, that is, the number of paths between any given pair of vertices is finite. Consider the category ${\rm rep}(Q)$ of locally finite dimensional $k$-representations of $Q$, which is a hereditary abelian $k$-category. For each $x\in Q_0$, one defines an indecomposable projective representation $P_x$ and an indecomposable injective representation $I_x$. Denote by ${\rm proj}(Q)$ and ${\rm inj}(Q)$ the additive full subcategories of ${\rm rep}(Q)$ generated by
the $P_x$ with $x\in Q_0$ and by the $I_x$ with $x\in Q_0$, respectively. As in the finite case, there exists a Nakayama equivalence $\nu: {\rm proj}(Q)\to {\rm inj}(Q)$. Furthermore, let ${\rm rep}^+\hspace{-1pt}(Q)$ be the full subcategory of $\rep(Q)$ generated by the representations which are presented by the projective representations in ${\rm proj}(Q)$, and let ${\rm rep}^-\hspace{-1pt}(Q)$ be the one generated by the representations which are co-presented by the injective representations in ${\rm inj}(Q).$ It is shown that ${\rm rep}^+\hspace{-1pt}(Q)$ and ${\rm rep}^-\hspace{-1pt}(Q)$ are Hom-finite, hereditary, abelian, and extension-closed in $\rep(Q)$; and that their intersection is the category 
of finite dimensional representations; see \cite[(1.15)]{BLP}.

\bigskip

Let $M$ be an indecomposable representation in $\rep(Q)$. If $M\in {\rm rep}^+\hspace{-1pt}(Q)$ has a minimal projective presentation
\vspace{-5pt}$$\xymatrixcolsep{14pt}\xymatrix{0\ar[r]& P_1 \ar[r]^f&
P_0 \ar[r] & M \ar[r]& 0,}$$ then we denote by $\DTr M$ the kernel of $\nu(f)$; and if $M\in {\rm rep}^-\hspace{-1pt}(Q)$ has a minimal injective
co-presentation \vspace{-8pt}$$\xymatrixcolsep{14pt}\xymatrix{0\ar[r]& M \ar[r] &
I_0 \ar[r]^g & I_1 \ar[r]& 0,}$$ then $\TrD M$ denotes the co-kernel of $\nu^-(g)\vspace{1pt}$. It is shown in \cite{BLP} that $\rep(Q)$ has an almost split sequence $\xymatrixcolsep{14pt}\xymatrix{0\ar[r] & \DTr M \ar[r] & N \ar[r] & M \ar[r] & 0}\hspace{-2pt}$ provided $M\in \rep^+\hspace{-1pt}(Q)$ is not projective; and an almost split sequence $\xymatrixcolsep{14pt}\xymatrix{0\ar[r] & M \ar[r] & N \ar[r] & \TrD M \ar[r] & 0}$
provided $M\in \rep^-\hspace{-1pt}(Q)$ is not injective. Now, we shall apply our previous results to give a short proof of the following result obtained in \cite[Section 3]{BLP}.

\bigskip

\begin{Theo} \label{rep+} Let $M$ be an indecomposable representation in $\rep(Q)$.

\begin{enumerate}[$(1)$]

\item If $M\in \rep^+\hspace{-1pt}(Q)$ is not projective, then
${\rm rep}^+\hspace{-1pt}(Q)$ has an almost split sequence $\xymatrixcolsep{14pt}\xymatrix{0\ar[r] & L \ar[r] & N \ar[r] & M \ar[r] & 0}\hspace{-2pt}$ if and only if $\,\DTr M$ is finite dimensional$\hspace{0.3pt};$  and in this case,  $L\cong \DTr M.$

\vspace{1pt}

\item If $M\in \rep^-\hspace{-1pt}(Q)$ is not injective, then
$\rep^-\hspace{-1pt}(Q)$ has an almost split sequence $\xymatrixcolsep{14pt}\xymatrix{0\ar[r] & M \ar[r] & N \ar[r] & L \ar[r] & 0}\hspace{-2pt}$ if and only if $\,\TrD M$ is finite dimensional$\hspace{0.3pt};$ and in this case, $L\cong \TrD M$.

\end{enumerate} \end{Theo}

\noindent{\it Proof.} We prove only Statement (1). Let $M\in \rep^+\hspace{-1pt}(Q)$ be not projective. Then,
$\rep(Q)$ has an almost split sequence $\delta: \xymatrixcolsep{15pt}\xymatrix{0\ar[r] & \DTr M \ar[r] & N \ar[r] & M \ar[r] & 0,}$ where $\DTr M\in \rep^-\hspace{-1pt}(Q)$; see \cite[(2.8)]{BLP}. If $\DTr M$ is finite dimensional, then $\delta$ lies in ${\rm rep}^+\hspace{-1pt}(Q)$, and hence it is an almost split sequence in ${\rm rep}^+\hspace{-1pt}(Q)$.

Conversely, assume that  $\xymatrixcolsep{15pt}\xymatrix{0\ar[r] & L \ar[r] & N \ar[r] & M \ar[r] & 0}\hspace{-2pt}$ is an almost split sequence in ${\rm rep}^+\hspace{-1pt}(Q)$. By Proposition \ref{ass-approx}, $\DTr M$ has a right injectively stable ${\rm rep}^+\hspace{-1pt}(Q)$-approximation $f: L\to \DTr M$. Suppose that $\DTr M$ is infinite dimensional. Then ${\rm supp}(\DTr M)$ contains a left infinite path; see \cite[(1.7)]{BLP}, which does not lie in ${\rm supp}\,L$; see \cite[(1.6)]{BLP}. In particular, there exists a vertex $x$ in ${\rm supp}(\DTr M)$ which does not lie in ${\rm supp}\,L$. Then $\Hom_{\rep(Q)}(P_x, L)=0\vspace{1pt}$, but there exists a non-zero morphism $g: P_x\to \DTr M$. Write $g=j h$, where $j: X\to \DTr M$ is a monomorphism and $h: P_x\to X$ is an epimorphism. Then, we have $\Hom_{\rep(Q)}(X, L)=0\vspace{0pt}$. Moreover, ${\rm supp}\,X$ is contained in the intersection of ${\rm supp}\,P_x$ and ${\rm supp}(\DTr M)$. Note that ${\rm supp}(\DTr M)$ has some vertices $a_1, \ldots, a_r$ such that every vertex in ${\rm supp}(\DTr M)$ is a predecessor of some of the $a_i$; see \cite[(1.6)]{BLP}. Since $Q$ is interval-finite, ${\rm supp}\,X$ is finite. Thus $X$ is finite dimensional. In particular, $X$ has an injective envelope $q: X\to J$ with $J\in {\rm inj}(Q)$. Observe that $\overline{j}$ factors through $\overline{\hspace{-1pt}f}$. Then $\bar{j}=\bar{0}$,
since $\Hom_{\rep(Q)}(X, L)=0\vspace{0pt}$. That is, $j$ is injectively trivial, and hence $j$ factors through $q$. On the other hand, since $\rep(Q)$ is hereditary and $\DTr M$ is indecomposable, $\Hom_{\rep(Q)}(J, \; \DTr M)=0.$ This yields $j=0$, a contradiction. The proof of the theorem is completed.

\bigskip\bigskip

{\sc Acknowledgement.} The first named author is supported in part by the Natu\-ral
Science and Engineering Research Council of Canada.

\end{document}